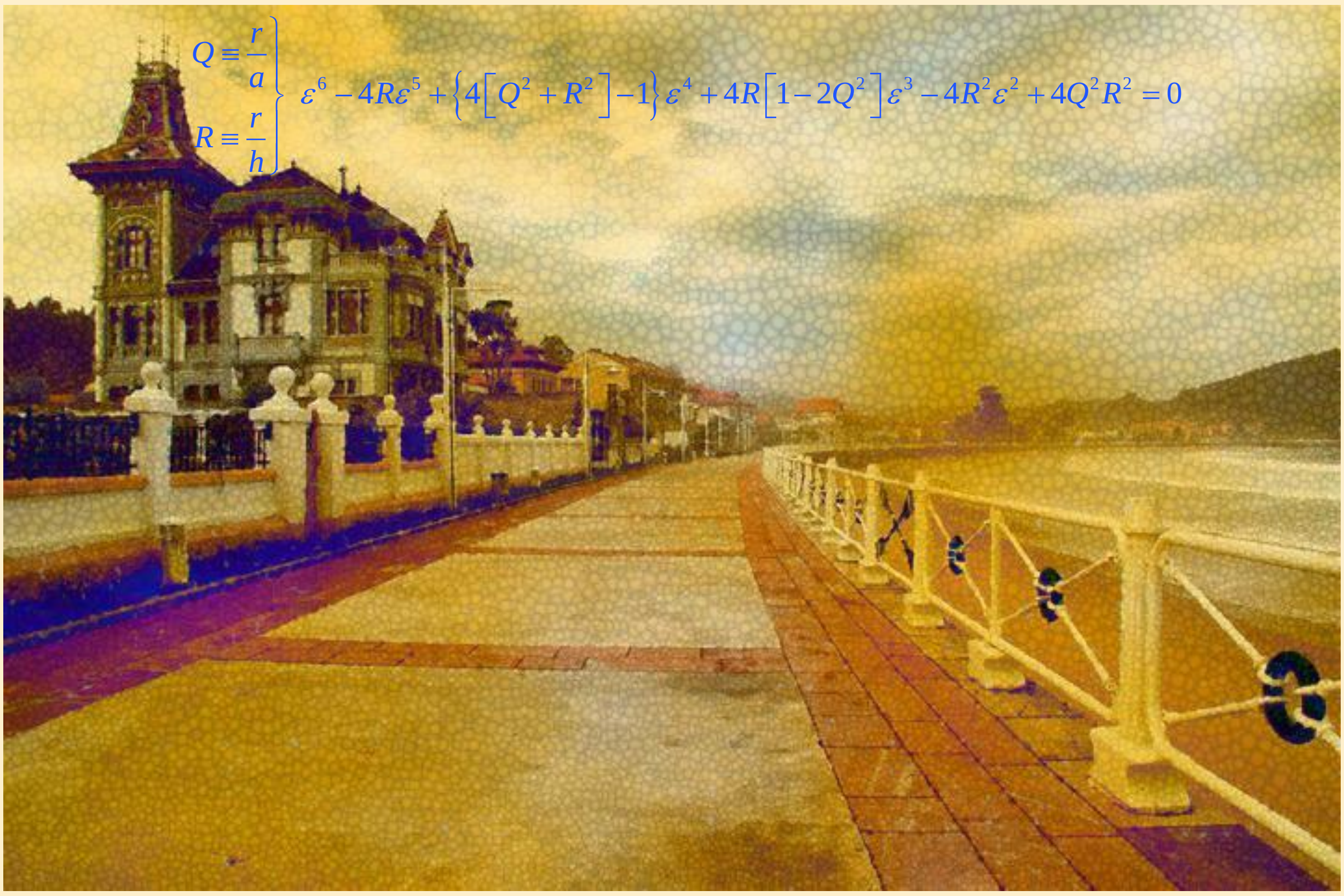

# Geometry Beyond Algebra

Hard – Hard = Easy

## The Theorem of Overlapped Polynomials (TOP) and its Application to the Sawa Masayoshi´s Sangaku Problem

The Adventure of Solving a Mathematical Challenge Stated in 1821

Jesús Álvarez Lobo
lobo@matematicas.net

Spain - 2011

In memoriam of Tere.

The solution presented here to the geometric Sangaku problem proposed by Masayoshi in 1821 and for which no solutions were known, don't belong to me, because I was inspired by an angel named Tere (Maria Teresa Suárez Fernández).

I offer this work as a tribute to the woman who decided to live beside me until the last day of her life, and now from Heaven she guides me with her Light and protects me with her Love.

I needed a way to tell her how much she means to me, until I can do it by myself, when at last she and I unite in Eternity.

Jesús Álvarez Lobo.

Oviedo, June 14, 2011.

ABSTRACT.

This work presents for the first time a solution to the 1821 unsolved Sawa Masayoshi´s problem, giving an *explicit and algebraically exact solution* for the *symmetric case* (particular case $b = c$, *i.e.*, $\triangle ABC \equiv$ right-angled isosceles triangle), see (1.60) and (1.61).

Despite the isosceles triangle restriction is not necessary, in view of the complexity of the explicit algebraic solution for the symmetric case, one can guessing the impossibility of achieving an explicit relationship for the *asymmetric case* (the more general case: $\triangle ABC \equiv$ right-angled scalene triangle). For this case is given a proof of *existence and uniqueness of solution* and a *proof of the impossibility* of getting such a relationship, even implicitly, if the sextic equation (2.54) it isn´t solvable. Nevertheless, in (2.56) – (2.58) it is shown the way to solve the asymmetric case under the condition that (2.54) be solvable.

Furthermore, it is proved that with a slight modification in the *final set of variables* ($\mathcal{F}$), it is still possible to establish a relation between them, see (2.59) and (2.61), which provides a bridge that connects the primitive relationship by means of *numerical methods*, for every given right-angled triangle $\triangle ABC$.

And as the attempt to solve Fermat's conjecture (or Fermat's last theorem), culminated more than three centuries later by Andrew Wiles, led to the development of powerful theories of more general scope, the attempt to solve the Masayoshi´s problem has led to the development of the *Theory of Overlapping Polynomials* (*TOP*), whose application to this problem reveals a great potential that might be extrapolated to other frameworks.

ACKNOWLEDGEMENTS

The first part of this work was presented at the University of Nagoya (Japan) on 5 August of 2011, by Hidetoshi Fukagawa, worldwide leading expert in Sangaku Geometry. It was also overseen by Tony Rothman (Princeton University), Kirsty Klafton (University of Edinburgh), Emanuele Delucchi (Bremen University, Germany), Peter Wong (Bates College, USA) and Miguel Amengual Covas (member of the commission of the European Mathematical Olympiad). I want to express to all of them my sincere appreciation.

I wish to express my gratitude in the most emphatic form to Alan Horwitz (Penn State University, USA) to assume the unpleasant task of checking my work, as an author authorized in arXiv to endorse me to submit my article.

**Introdution.**

Proposed by the Nippon mathematician *Sawa Masayoshi* in 1821, **this problem remained unsolved.**

As shown in the picture below, an ellipse $\mathcal{E}$ is inscribed in $\triangle ABC$, right-angled triangle at $A$, with its major axis parallel to the hypotenuse $a$ of $\triangle ABC$. Within the ellipse are inscribed two circumferences of radius $r$, $\mathcal{C}_2$ and $\mathcal{C}_3$, that are tangent to the ellipse and each other at the center of the ellipse. Inside-tangent to the sides $b$ and $c$ of the right angle $A$ and outside-tangent to the ellipse, there is another circumference, $\mathcal{C}_1$, with the same radius $r$.
Challenge: find a relationship between $a$, $b$, $c$ and $r$.

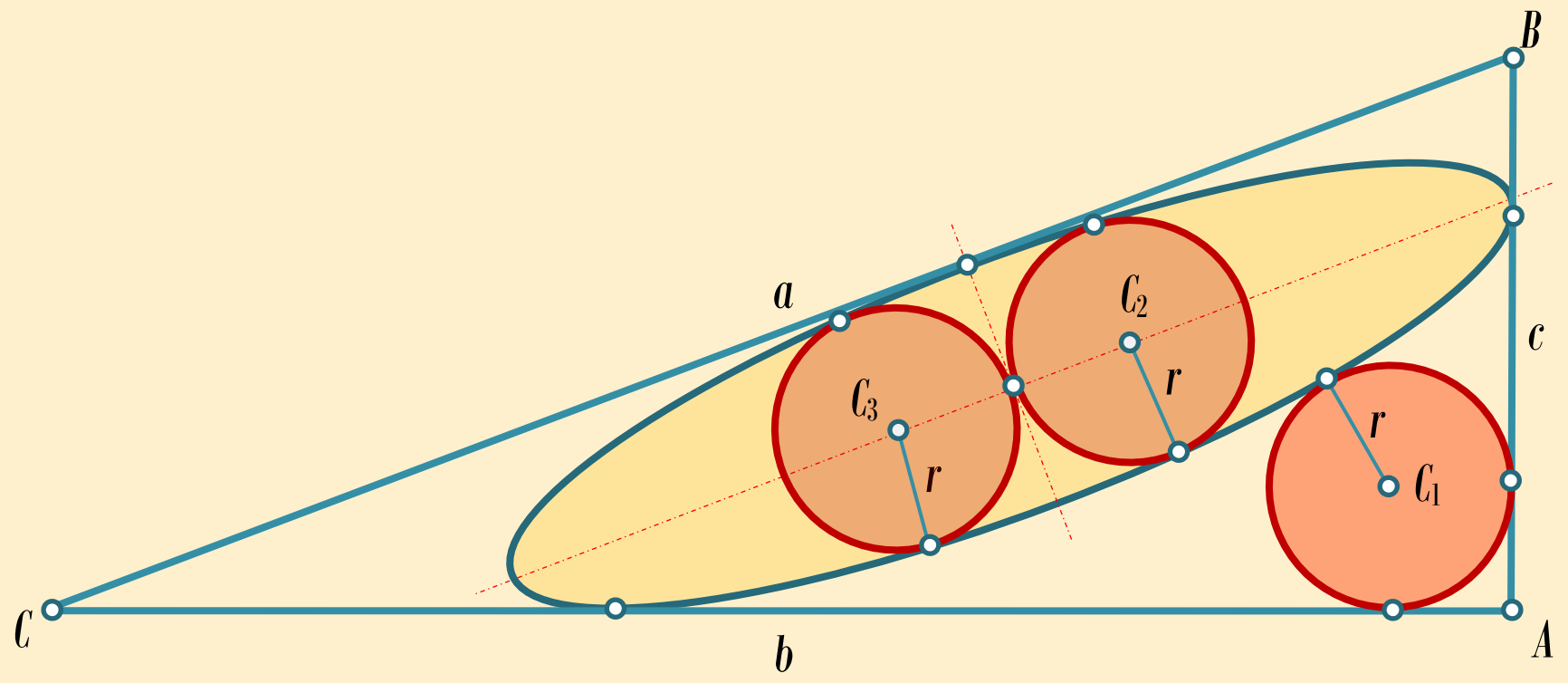

The Travel Diary of Mathematician Yamaguchi Kanzan. "I have arrived at Tatsuno, city near Himeji, and visited the Syosya temple to record a sangaku proposed by Sawa in 1821. In the evening, visited me and he showed me an unsolved problem and two sangaku problems of the Syosya temple. I have written them down in my diary as follows: . . ."*From the book "Sacred Mathematics. Japanese Temple Geometry"*, by Fukagawa Hidetoshi and Tony Rothman.

According to the brilliant mathematician Paul Erdös, "*A math problem that takes more than a century to be solved, it is a number theory problem*". And I would add: "*If a seemingly simple mathematical problem, because it hides an underlying extreme complexity, takes more than a century to be solved, most likely it is that it will be a number theory problem... or a geometric Sangaku problem, and only an exhaustive analysis could lead to its solution*". Jesús Álvarez Lobo.

The Masayoshi´s problem kept a surprise: the explicit relationship between the length of one side of $\triangle ABC$ and the radius of the circumferences or is so extensive and intricate that it cannot be written without fragmentation, or both the path to prove the impossibility of getting a relationship for the asymmetric case and the equations for the alternative case are huge and of an bewildering complexity, despite the valuable aid provided by the TOP. Everybody can believe in the existence of more simply relationships, (that´s easy!), but the reality is that, up today, the ones presented here are the unique known.

Some approaches, cut-off paths and other collateral interesting results found through the surveying, has been moved to final appendixes, in order to lighten the reading.

In spite of the proofs are carefully developed, some steps have been omitted for brevity. However, it isn´t hard to follow the exhaustive evolving of the work.

## §0. Theory of Overlapped Polynomials.

### 0.1. Definitions.

**Definition 1.** *Class* $\mathcal{O}_n(x-r)^{\mu}$ *of Overlapped Polynomials: Set of all polynomials of degree n in the variable x which have the root x = r with multiplicity not lower than* $\mu$.

**Definition 2.** *Monic Transformed of a Polynomial:*

$$\mathrm{M}\left(\sum_{i=0}^{n} a_i x^i\right) = \frac{1}{a_n}\sum_{i=0}^{n} a_i x^i \qquad (0.01)$$

### 0.2. *Lemma of Overlapped Polynomials 1 (LOP1).*

*The difference of the monic transforms of two overlapped polynomials of degree n in the variable x, that share the root r with multiplicity not lower than* $\mu$*, it is a polynomial of degree m < n that keeps the root r with multiplicity at least* $\mu$:

$$S(x),\ T(x) \in \mathcal{O}_n(x-r)^{\mu} \Rightarrow \mathrm{M}\{S(x)\} - \mathrm{M}\{T(x)\} \in \mathcal{O}_m(x-r)^{\mu},\ m<n \qquad (0.02)$$

**Proof.**

$$\left.\begin{array}{l} \mathrm{M}\{S(x)\} = (x-r)^{\mu} s(x) \\ \mathrm{M}\{T(x)\} = (x-r)^{\mu} t(x) \end{array}\right\} \Rightarrow \mathrm{M}\{S(x)\} - \mathrm{M}\{T(x)\} = (x-r)^{\mu}\left[s(x)-t(x)\right] \qquad (0.03)$$

where $s$ and $t$ are monic polynomials of degree $n-\mu$,

$$\begin{array}{l} s(x) = x^{n-\mu} + \displaystyle\sum_{i=1}^{n-\mu} a_{n-\mu-i} x^{n-\mu-i} \\ t(x) = x^{n-\mu} + \displaystyle\sum_{i=1}^{n-\mu} b_{n-\mu-i} x^{n-\mu-i} \end{array} \qquad (0.04)$$

$$\therefore\ s(x)-t(x) = \sum_{i=1}^{n-\mu} (a_{n-\mu-i} - b_{n-\mu-i}) x^{n-\mu-i} \qquad (0.05)$$

Therefore,

$$\mathrm{M}\{S(x)\} - \mathrm{M}\{T(x)\} = (x-r)^{\mu} \sum_{i=1}^{n-\mu} (a_{n-\mu-i} - b_{n-\mu-i}) x^{n-\mu-i} \in \mathcal{O}_{m<n}(x-r)^{\mu} \qquad (0.06)$$

*Be noted that, by Viète-Girard formula, the degree of S(x) - T(x) is less than n - 1 iff the sum of the roots of s(x) is equal to the sum of the roots of t(x), the degree of S(x) - T(x) is less than n - 2 iff also the sum of the binary products of the roots of s(x) is equal to the sum of the binary products of the roots of t(x), and so on...* (0.07)

### 0.3. *Lemma of Overlapped Polynomials 2 (LOP2).*

$$S,\ T \in \mathcal{O}_n(x-r)^{\mu} \Rightarrow \Delta_{ST}^{\mu-1} \equiv \mathrm{M}\left\{S^{(\mu-1}\right\} - \mathrm{M}\left\{T^{(\mu-1}\right\} \in \mathcal{O}_k(x-r)^1,\ k \le n-\mu \qquad (0.08)$$

denoting $S^{(\mu-1}$, $T^{(\mu-1}$ the derivatives of order $\mu$ - 1 of *S* and *T*, respectively.

**Proof.**

*S*(*x*) and *T*(*x*) can be factored as in (0.03), $S(x) = (x\text{-}r)^{\mu} s(x)$, $T(x) = (x\text{-}r)^{\mu} t(x)$, where *s* and *t* are the monic polynomials defined in (0.03).

Applying the *Leibniz´s formula* for the $\mu - 1$ - th derivative of the product to *S*(*x*),

$$\{S(x)\}^{(\mu-1} = \{(x-r)^{\mu} s\}^{(\mu-1} = \sum_{i=0}^{\mu-1} \binom{\mu-1}{i} \{(x-r)^{\mu}\}^{(\mu-1-i} s^{(i} \qquad (0.09)$$

and similarly for *T*(*x*),

$$\{T(x)\}^{(\mu-1} = \{(x-r)^{\mu} t^{(\mu-1}\} = \sum_{i=0}^{\mu-1} \binom{\mu-1}{i} \{(x-r)^{\mu}\}^{(\mu-1-i} t^{(i} \qquad (0.10)$$

It is straightforward that the terms of the sums (0.10) and (0.11) are all of degree $n - \mu + 1$ and contain the powers of $x - r$ with integer exponents between 1 and $\mu$, so,

$$\{S(x)\}^{(\mu-1} = (x-r)p(x) \qquad \{T(x)\}^{(\mu-1} = (x-r)q(x) \qquad (0.11)$$

where *p* and *q* are polynomials of degree $n - \mu$,

$$p(x) = \sum_{k=0}^{n-\mu} p_k x^k \qquad q(x) = \sum_{k=0}^{n-\mu} q_k x^k \qquad (0.12)$$

Hence,

$$\mathrm{M}\left\{S^{(\mu-1}\right\} = \mathrm{M}\left\{(x-r)p(x)\right\} = (x-r)\mathrm{M}\left\{p(x)\right\} = \frac{1}{p_{n-\mu}}(x-r)\sum_{i=0}^{n-\mu} p_i x^i \qquad (0.13)$$

$$\mathrm{M}\left\{T^{(\mu-1}\right\} = \mathrm{M}\left\{(x-r)q(x)\right\} = (x-r)\mathrm{M}\left\{q(x)\right\} = \frac{1}{q_{n-\mu}}(x-r)\sum_{i=0}^{n-\mu} q_i x^i \qquad (0.14)$$

$$\mathrm{M}\left\{S^{(\mu-1}\right\} - \mathrm{M}\left\{T^{(\mu-1}\right\} = (x-r)\left(\frac{1}{p_{n-\mu}}\sum_{i=0}^{n-\mu} p_i x^i - \frac{1}{q_{n-\mu}}\sum_{i=0}^{n-\mu} q_i x^i\right) \qquad (0.15)$$

But the summations of (0.15) are monic polynomials of degree *n* - $\mu$, therefore, applying the LOP1,

$$\Delta_{ST}^{\mu-1} \equiv \mathrm{M}\left\{S^{(\mu-1}\right\} - \mathrm{M}\left\{T^{(\mu-1}\right\} = (x-r)\left(\frac{1}{p_{n-\mu}}\sum_{i=0}^{n-\mu-1} p_i x^i - \frac{1}{q_{n-\mu}}\sum_{i=0}^{n-\mu-1} q_i x^i\right) \in \mathcal{O}_k(x-r)^1,\ k \le n-\mu \qquad (0.16)$$

## 0.4. *Theorem of Overlapped Polynomials (TOP).*

This theorem is an extension of LOP2, removing the restriction of equal degree of the overlapped polynomials involved. Therefore, the TOP represents the widest generalization of the Theory of Overlapped Polynomials.

The theorem is stated as follows:

$$\left.\begin{array}{l} S \in \mathcal{O}_n(x-r)^{\sigma} \\ T \in \mathcal{O}_m(x-r)^{\tau} \\ \delta \equiv n-m \in \mathbb{N} \\ \mu \equiv \min\{\sigma, \tau\} \end{array}\right\} \Rightarrow \Delta_{ST}^{\mu-1} \equiv \mathrm{M}\left\{S^{(\mu-1)}\right\} - \mathrm{M}\left\{(x^{\delta}T)^{(\mu-1)}\right\} \in \mathcal{O}_k(x-r)^1,\ k \leq n-\mu \tag{0.17}$$

**Proof.**

$$\left.\begin{array}{l} S \in \mathcal{O}_n(x-r)^{\sigma} \\ T \in \mathcal{O}_m(x-r)^{\tau} \\ \delta \equiv n-m \in \mathbb{N} \\ \mu \equiv \min\{\sigma, \tau\} \end{array}\right\} \Rightarrow \left\{\begin{array}{l} S \in \mathcal{O}_n(x-r)^{\mu} \subseteq \mathcal{O}_n(x-r)^{\sigma} \\ x^{\delta}T \in \mathcal{O}_{m+\delta}(x-r)^{\mu} = \mathcal{O}_n(x-r)^{\mu} \subseteq \mathcal{O}_n(x-r)^{\tau} \end{array}\right. \tag{0.18}$$

where the inclusion relationship are due to the *nested* characterization stated in the definition of the *overlapped classes of polynomials*. Hence, applying the LOP2,

$$S,\ x^{\delta}T \in \mathcal{O}_n(x-r)^{\mu} \Rightarrow \Delta_{ST}^{\mu-1} \equiv \mathrm{M}\left\{S^{(\mu-1)}\right\} - \mathrm{M}\left\{(x^{\delta}T)^{(\mu-1)}\right\} \in \mathcal{O}_k(x-r)^1,\ k \leq n-\mu \tag{0.19}$$

## 0.4. Applications of the TOP.

In the Masayoshi´s problem, to relate the radius of the circumferences inscribed in the ellipse with the radius of the outside tangent circle to the ellipse, is necessary to calculate the coordinates of the tangent point *T* as a function only of *b*, *c* and *r*, but the equations that determine it are complete quartics of a extreme complexity. On the other hand, as should be shown below, it is possible to define three independent quartics in the abscise $x_T$ of *T*. At first glance, it could look like superfluous to state three different condition for the same purpose, but amazingly it is possible to get advantage of this circumstance, what can be crucial in certain conditions.

In general, there may be several possibilities, depending on the degree of multiplicity of the root shared and on the depth (order of the derivative) at which the theorem be applied.

The tangent point *T* of two curves is twofold and this reasoning led to the development of the *Theory of Overlapping Polynomials* in the generalized form in which it is presented, i.e., considering *multiplicities* of the roots. However, we must consider that multiplicity cannot be inferred from the mere consideration of the geometry of the problem (Algebra is generous and often gives more than is asked!), since in the *transformation of rationalization* (squaring) *strange solutions* can be introduced. At this problem the TOP has been implemented in its simplest form (considering *simple roots*), but even so a unit reduction in the degree of an equation can be crucial to achieve a solution, and will always represent a substantial simplification of the calculations and the solution.

Moreover, the TOP can be applied repeatedly, thereby achieving further reductions in the degree of the equations. For instance, from three polynomial equations of any degree bigger than 2 that share a root $x = r$, by successive applications of the TOP, is possible to reduce the degree of the polynomial equation of least degree in two units.

For example, suppose you have two equations of 6th degree and one of 5th degree for the determination of a certain magnitude $s$. At first, the problem should be algebraically unsolvable, but it might go beyond the limits of Algebra by applying the TOP: with a first application, it is we obtained a *quintic* from the two *sextics*, and then, applied to the two *quintics* gives a *quartic* which is already solvable algebraically.

Suppose now that you have three equations of 6th grade for the determination of a certain magnitude $s$. Since the sextics are, in general, algebraically unsolvable, it looks like that the problem hasn´t algebraic solution, but it possible trespass the limits of Algebra by applying the TOP: applied at a pair of sextics (for instance, to the first a second) we get a quintic, and applied to other pair of sextics (for example, first and third) we get another quintic. Finally, the application of the TOP to the pair of quintics derived from the sextics gives us a quartic that is always solvable by radicals.

It is not particularly unusual to find different equations for the determination of a value of a variable under certain conditions, as problems often lets different approaches leading to different equations. Usually the equation will be of the same degree as they have to contain all the necessary information to define the problem, but even if they have different degrees it is possible to apply the TOP.

Let´s see a concrete example to clarify the way the TOP is applied.

Suppose you have three polynomials,

$$u(x) \equiv 9x^6 - 24x^5 + 40x^4 - 38x^3 + 33x^2 - 14x + 2 \tag{0.20}$$

$$v(x) \equiv 3x^5 + 11x^4 - 70x^3 + 10x^2 + 67x - 21 \tag{0.21}$$

$$w(x) \equiv 6x^5 - 17x^4 + 14x^3 - 4x + 1 \tag{0.22}$$

that share certain unknown root, $x = r$.

Applying the TOP to $u$ and $v$,

$$\left.\begin{array}{l} u \in \mathcal{O}_6(x-r) \\ v \in \mathcal{O}_5(x-r) \\ \delta \equiv 6-5=1 \\ \mu \equiv 1 \end{array}\right\} \Rightarrow \Delta^0_{uv} \equiv \mathrm{M}\{u\} - \mathrm{M}\{xv\} \in \mathcal{O}_k(x-r),\ k \le 6-1=5 \tag{0.23}$$

i.e.,

$$\Delta^0_{u(xw)} = \frac{1}{9}u - \frac{1}{6}xv = \frac{1}{6}x^5 + \frac{19}{9}x^4 - \frac{38}{9}x^3 + \frac{13}{3}x^2 - \frac{31}{18}x + \frac{2}{9} \tag{0.24}$$

Now we have three polynomials of 5th degree that share the root $x = r$, and applying the TOP to the three pairs of quintic equations we will obtain three quartic equations that preserves the same common root.

Let us denote with $u \bullet w$ the quintic inherited by means of the TOP from the sextic $u$ and the quintic $w$.

$$(u\bullet w)\bullet v \equiv \mathrm{M}\{u\bullet w\}-\mathrm{M}\{v\}=6u\bullet w-\frac{1}{3}v=9x^4-2x^3+\frac{68}{3}x^2-\frac{98}{3}x+\frac{25}{3} \quad (0.25)$$

$$(u\bullet w)\bullet w \equiv \mathrm{M}\{u\bullet w\}-\mathrm{M}\{w\}=6u\bullet w-\frac{1}{6}w=\frac{31}{2}x^4-\frac{83}{3}x^3+26x^2-\frac{29}{3}x+\frac{7}{6} \quad (0.26)$$

$$v\bullet w \equiv \Delta^0_{vw}=\mathrm{M}\{v\}-\mathrm{M}\{w\}=\frac{1}{3}v-\frac{1}{6}w=\frac{13}{2}x^4-\frac{77}{3}x^3+\frac{10}{3}x^2+23x-\frac{43}{6} \quad (0.27)$$

Now we have three polynomials of 4th degree that share the root $x = r$, and applying the TOP to the three pairs of quartic equations we will obtain three cubic equations that preserves the root $x = r$,

$$[(u\bullet w)\bullet v][(u\bullet w)\bullet w] \equiv \mathrm{M}\{(u\bullet w)\bullet v\}-\mathrm{M}\{(u\bullet w)\bullet w\}=\frac{436}{279}x^3+\frac{704}{837}x^2-\frac{2516}{837}x+\frac{712}{837} \quad (0.28)$$

$$[(u\bullet w)\bullet v][v\bullet w]=\mathrm{M}\{(u\bullet w)\bullet v\}-\mathrm{M}\{v\bullet w\}=\frac{436}{117}x^3+\frac{704}{351}x^2-\frac{2516}{351}x+\frac{712}{351} \quad (0.29)$$

$$[(u\bullet w)\bullet w][v\bullet w]=\mathrm{M}\{(u\bullet w)\bullet w\}-\mathrm{M}\{v\bullet w\}=\frac{872}{403}x^3+\frac{1408}{1209}x^2-\frac{5032}{1209}x+\frac{1424}{1209} \quad (0.30)$$

If we would try to implement now the TOP to any pair of the above cubic equations, for instance,

$$\mathrm{M}\{[(u\bullet w)\bullet v][(u\bullet w)\bullet w]\}-\mathrm{M}\{[(u\bullet w)\bullet v][v\bullet w]\}\equiv 0 \quad (0.31)$$

we would get the null polynomial because the three equations are equivalent. So, it is impossible to further reduce the degree of the equations.

Solving any of the cubic equations, for instance, get,

$$x=\pm\frac{\sqrt{86633}}{218}-\frac{95}{218},\quad x=\frac{1}{3} \quad (0.32)$$

one of whose roots is that which is shared by the three initial equations.

To determine which of the three roots is what is seeking is necessary to use some criteria inferred from the conditions of the problem. In this case, the root shared is $x = 1/3$.

## §1. Symmetric case: $b = c$.

In the orthonormal Cartesian reference system provided in the form shown in the figure below, the equation of the ellipse is reduced to its canonical form and geometric configuration has axial symmetry on the coordinate axes.

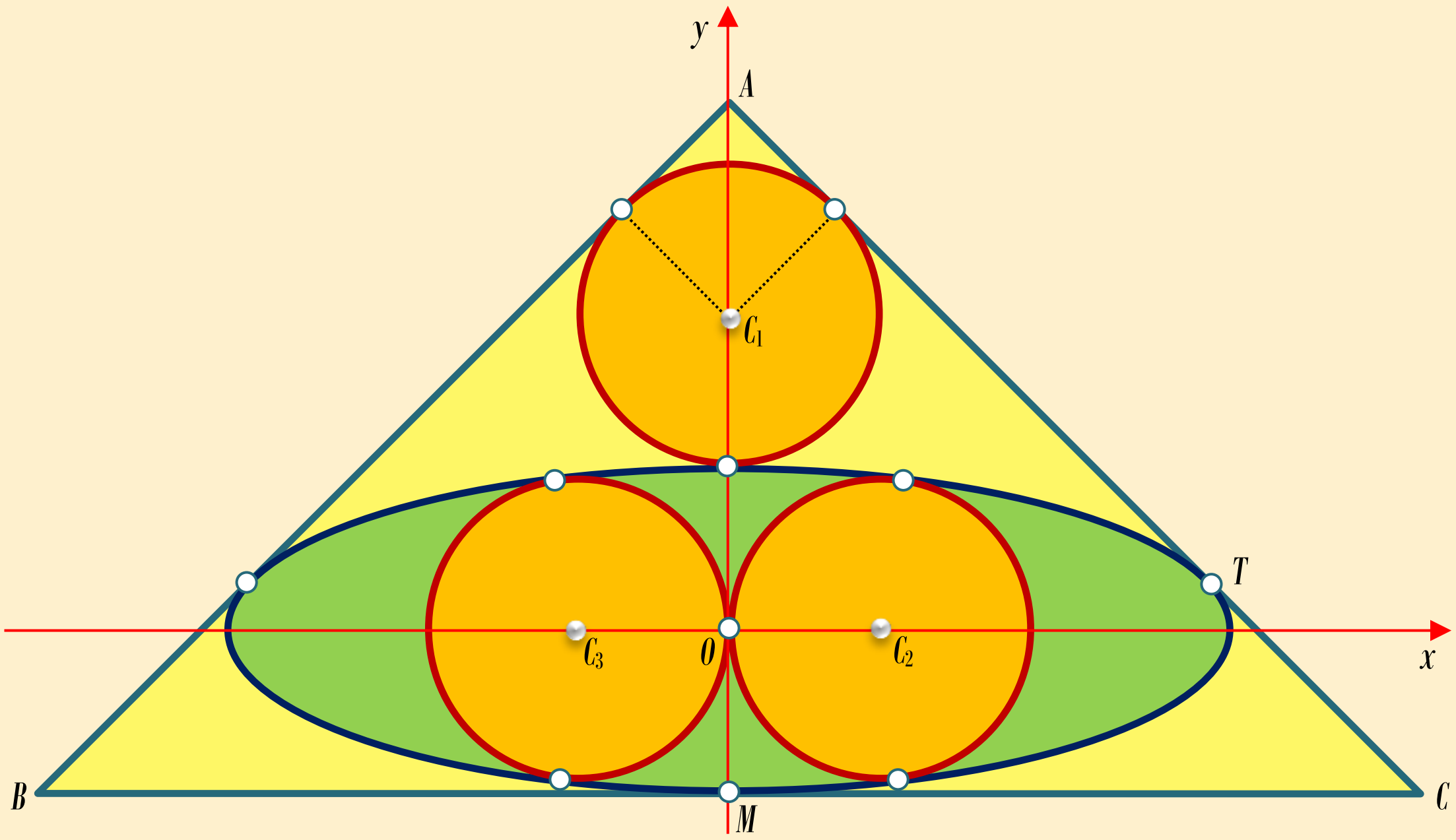

**Equation of the ellipse in canonical form** (center at the origin and major-axis of the ellipse coincident with the *x*-axis):

$$\mathscr{E} \equiv \frac{x^2}{\alpha^2} + \frac{y^2}{\beta^2} = 1 \qquad (1.01)$$

Let $t$ be the straight line defined by $A$ and $C$:

$$t \equiv y = -x + y_A \qquad (1.02)$$

From the intersection of $t$ and the ellipse $\mathscr{E}$:

$$t \cap \mathscr{E} \equiv \frac{x^2}{\alpha^2} + \frac{(-x+y_A)^2}{\beta^2} = 1,$$

that can be written as,

$$\left(1 + \frac{\beta^2}{\alpha^2}\right)x^2 - 2y_A x + y_A^2 - \beta^2 = 0 \qquad (1.03)$$

However, the necessary and sufficient condition for $t$ to be tangent to $\mathscr{E}$ at $T$ is that the last equation has twofold roots, or, equivalently, discriminant null:

$$t \cap \mathscr{E} = T \Leftrightarrow \Delta = (2y_A)^2 - 4\left(1 + \frac{\beta^2}{\alpha^2}\right)(y_A^2 - \beta^2) = 0,$$

or,

$$y_A^2 - \left(1 + \frac{\beta^2}{\alpha^2}\right)(y_A^2 - \beta^2) = 0 \Leftrightarrow -\frac{\beta^2}{\alpha^2}y_A^2 + \left(1 + \frac{\beta^2}{\alpha^2}\right)\beta^2 = 0 \Leftrightarrow y_A = \pm\sqrt{\alpha^2 + \beta^2},$$

and for the positive semi-ellipse,

$$y_A = +\sqrt{\alpha^2+\beta^2} \tag{1.04}$$

Furthermore, straight by simple observation of the picture are obtained the coordinates of *A*,

$$A(x_A, y_A) = A\left(0, \beta + r + r\sqrt{2}\right)$$

as $\overline{AC_1}$ is the diagonal of a square of side *r*. Hence,

$$y_A = (1+\sqrt{2})r + \beta \tag{1.05}$$

The pair of equations (1.04) and (1.05), that arises in a fortunate choice of the reference system, will be crucial to the resolution of the problem. Equating these two expressions of $y_A$, we have:

$$\sqrt{\alpha^2+\beta^2} = \beta + (1+\sqrt{2})r \Leftrightarrow \alpha^2 + \cancel{\beta^2} = \cancel{\beta^2} + 2\beta(1+\sqrt{2})r + (1+\sqrt{2})^2 r^2,$$

that can be expressed as

$$\frac{1}{1+\sqrt{2}}\alpha^2 = 2\beta r + (1+\sqrt{2})r^2 \tag{1.06}$$

In order to compact the expression, let us denote

$$p \equiv \frac{1}{1+\sqrt{2}} = \sqrt{2}-1 \tag{1.07}$$

Then, (1.06) is written as,

$$p\alpha^2 = 2\beta r + \frac{r^2}{p}, \tag{1.08}$$

and from this one,

$$p^2 = 2p\frac{\beta}{\alpha^2}r + \frac{r^2}{\alpha^2} \tag{1.09}$$

**Intersection of $\mathcal{C}_2$ and $\mathcal{E}$:**

$$\left.\begin{aligned} \mathcal{E} &\equiv \frac{x^2}{\alpha^2} + \frac{y^2}{\beta^2} = 1 \\ \mathcal{C}_2 &\equiv (x-r)^2 + y^2 = r^2 \end{aligned}\right\} \mathcal{E} \cap \mathcal{C}_2 \equiv \frac{x^2}{\alpha^2} + \frac{r^2-(x-r)^2}{\beta^2} = 1,$$

quadratic equation in *x* that can be written as,

$$\left(1-\frac{\beta^2}{\alpha^2}\right)x^2 - 2rx + \beta^2 = 0 \tag{1.10}$$

Now, taking into account that, by definition, the eccentricity of the ellipse is

$$\varepsilon = \frac{\gamma}{\alpha} = \frac{\sqrt{\alpha^2 - \beta^2}}{\alpha} = \sqrt{1 - \frac{\beta^2}{\alpha^2}} \tag{1.11}$$

the equation (1.10) can be expressed as,

$$\varepsilon^2 x^2 - 2rx + \beta^2 = 0 \tag{1.12}$$

Moreover, the **condition of tangency** between the circle with center at $C_2$ and the ellipse is that the contact occurs only between two points of equal abscissa, what will happen iff the discriminant of (1.12), as a quadratic equation in $x$, is null, and this condition is translated in the following relationship:

$$\Delta = 0 \Leftrightarrow (2r)^2 - 4\varepsilon^2\beta^2 = 0,$$

that implies,

$$r = \beta\varepsilon \tag{1.13}$$

However, should be noted that the maximum radius of the circles inscribed in the semi-ellipse of semi-major axis $\alpha$ is $\alpha/2$. Whereas as a constant, $r$ can be expressed as function of $\varepsilon$ as the only variable:

$$\left.\begin{aligned} & r = \beta\varepsilon \\ & 1 - \frac{\beta^2}{\alpha^2} = \varepsilon^2 \Leftrightarrow \beta^2 = \alpha^2(1-\varepsilon^2) \end{aligned}\right\} \Rightarrow r = \pm\alpha\varepsilon\sqrt{1-\varepsilon^2},$$

but in the context of the problem,

$$r(\varepsilon) = \alpha\varepsilon\sqrt{1-\varepsilon^2} \tag{1.14}$$

Let us study the increasing / decreasing of this function: $r'(\varepsilon) = \alpha\dfrac{1-2\varepsilon^2}{\sqrt{1-\varepsilon^2}} = 0 \Rightarrow \varepsilon = \pm\dfrac{1}{\sqrt{2}}$,

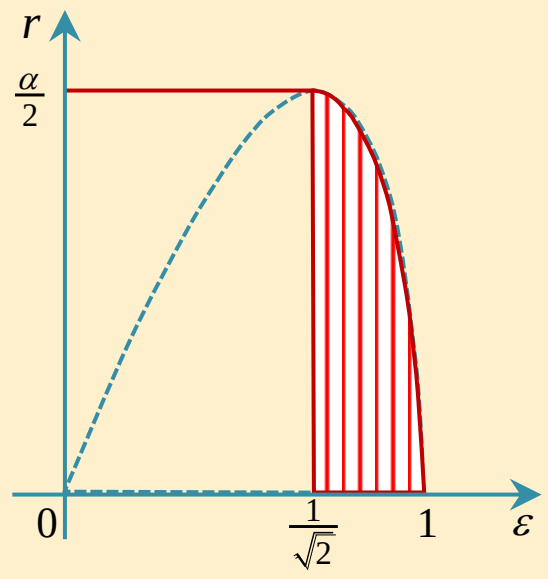

$$\begin{cases} r'(\varepsilon) > 0, \; 0 < \varepsilon < \dfrac{1}{\sqrt{2}} \\ r'(\varepsilon) < 0, \; \dfrac{1}{\sqrt{2}} < \varepsilon < 1 \end{cases} \tag{1.15}$$

as can be checked trivially.

So, $r$ reaches a *relative and absolute maximum* for $\varepsilon = \dfrac{1}{\sqrt{2}}$, of value $r\left(1/\sqrt{2}\right) = \alpha\dfrac{1}{\sqrt{2}}\sqrt{1-\left(\dfrac{1}{\sqrt{2}}\right)^2} = \dfrac{\alpha}{2}$, i.e.,

$$r(\varepsilon)_{\max} = r\left(1/\sqrt{2}\right) = \frac{\alpha}{2} \tag{1.16}$$

Hence, the domain of $r(\varepsilon)$ is,

$$Dom\, r(\varepsilon) = \left(\frac{1}{\sqrt{2}}, 1\right) \tag{1.17}$$

The left end of the domain of $r(\varepsilon)$ coincides with the minimum value of the eccentricity for which the radius of the circles inscribed in the semi-ellipse reaches a maximum. This limit value can also be determined as follows:

$$r = \beta\varepsilon = \frac{\alpha}{2} \Leftrightarrow \frac{\beta^2}{\alpha^2} = \frac{1}{4\varepsilon^2} \Leftrightarrow 1-\varepsilon^2 = \frac{1}{4\varepsilon^2},$$

i.e.,

$$\boxed{4\varepsilon^4 - 4\varepsilon^2 + 1 = 0} \qquad (1.18)$$

biquadratic equation whose unique real and positive root is $\varepsilon = 1/\sqrt{2}$.

The circle with center at $C_2$ inscribed in the right semi-ellipse has maximum radius $r = \alpha/2$ for $\varepsilon \le 1/\sqrt{2}$, as it only has contact with the ellipse at the right apex (point of tangency). If $\varepsilon = 1$ the ellipse is degenerated in a segment.

In short,

$$\boxed{r = \begin{cases} \beta\varepsilon \text{ if } \varepsilon \in \left(\frac{1}{\sqrt{2}},\, 1\right) \\ \frac{\alpha}{2} \text{ if } \varepsilon \in \left[0,\, \frac{1}{\sqrt{2}}\right] \end{cases}} \qquad (1.19)$$

But $r = \alpha/2$ cannot be a solution of the problem.

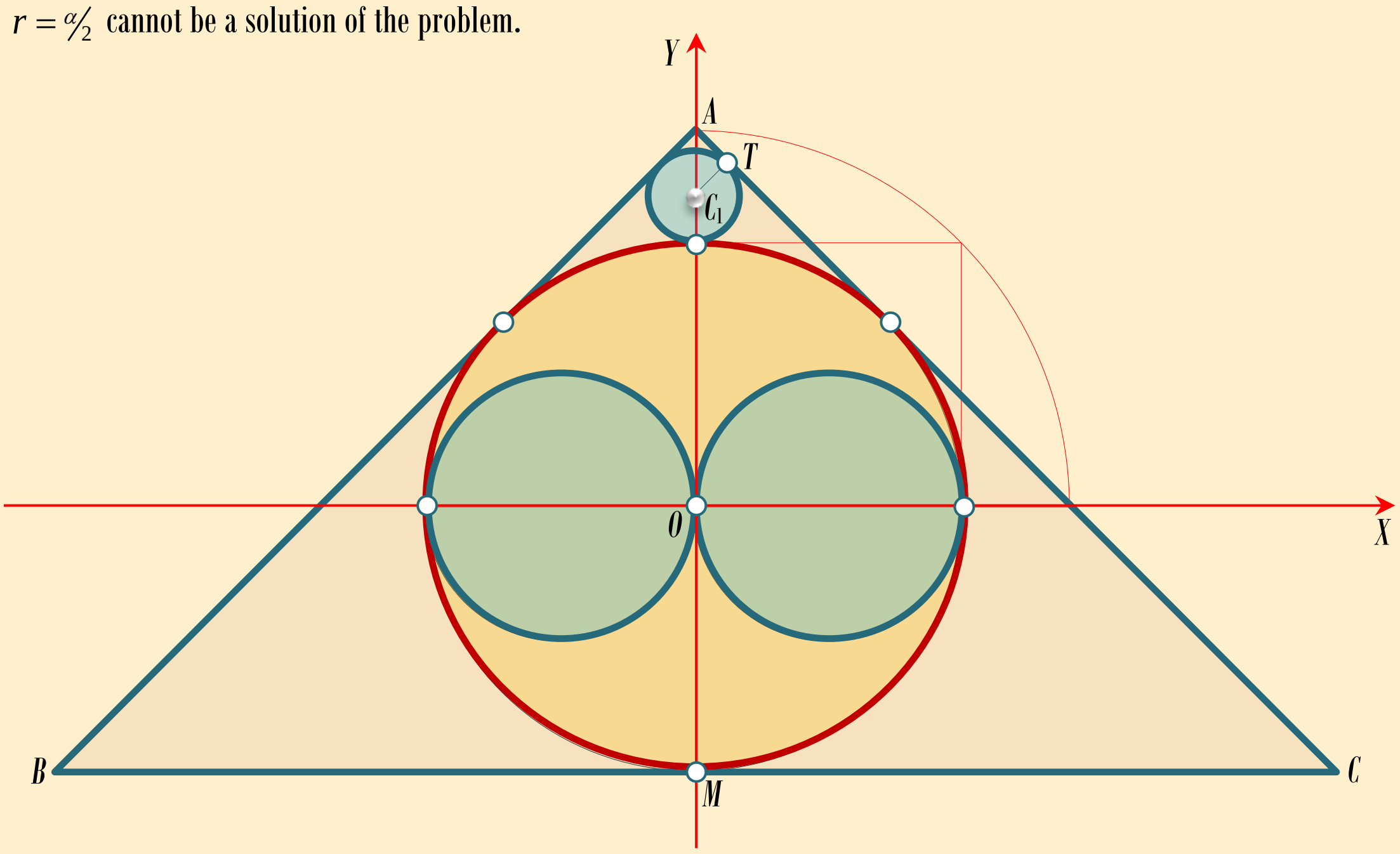

Indeed, referring to the figure above, according to (1.04) and (1.05), $y_A = \overline{OA} = \sqrt{\alpha^2 + \beta^2} = (1+\sqrt{2})r_1 + \beta$. Hence,

$$\alpha^2 + \cancel{\beta^2} = (1+\sqrt{2})^2 r_1^2 + 2(1+\sqrt{2})r_1 + \cancel{\beta^2},$$

$$\boxed{(1+\sqrt{2})^2 r_1^2 + 2(1+\sqrt{2})r_1 - \alpha^2 = 0} \qquad (1.20)$$

quadratic equation in $r_1$, whose solutions are,

$$r_1 = \pm\left(\sqrt{2}-1\right)\left(\sqrt{\alpha^2+1} \mp 1\right).$$

Discarding the negative root, that has no sense in this context,

$$r_1 = \left(\sqrt{2}-1\right)\left(\sqrt{\alpha^2+1}-1\right) \tag{1.21}$$

and from (1.21), we obtain:

$$r_1 < \frac{1}{2}\left(\sqrt{\alpha^2+1}-1\right) < \frac{1}{2}\left(\sqrt{\alpha^2+2\alpha+1}-1\right) = \frac{1}{2}\left(\sqrt{(\alpha+1)^2}-1\right) = \frac{1}{2}\left((\alpha+1)-1\right) = \frac{\alpha}{2} = r_2 = r_3.$$

Consequently, the solution of the problem, unique due the monotony of $r(\varepsilon)$, is confined to the domain of $r(\varepsilon)$, established in (1.17).

Now, eliminating $r$ in (1.09) through the relationship (1.13), which links the radius of the inside the ellipse circunferences with the semi-major axis of the ellipse and its eccentricity, we obtain:

$$\left.\begin{array}{l} p^2 = 2p\dfrac{\beta}{\alpha^2}r + \dfrac{r^2}{\alpha^2} \\ r = \beta\varepsilon \end{array}\right\} \Rightarrow p^2 = 2p\frac{\beta^2}{\alpha^2}\varepsilon + \frac{\beta^2}{\alpha^2}\varepsilon^2 \Leftrightarrow p^2 = 2p(1-\varepsilon^2)\varepsilon + (1-\varepsilon^2)\varepsilon^2,$$

from which it is obtained the complete quartic equation in $\varepsilon$,

$$\varepsilon^4 + 2p\varepsilon^3 - \varepsilon^2 - 2p\varepsilon + p^2 = 0 \tag{1.22}$$

By the lineal *Tschirnhausen transformation*,

$$\varepsilon \equiv \xi - p/2 \tag{1.23}$$

is eliminated the term of third grade in (1.22) and, after replacing $p$ by its value (1.07) and simplifying, we get the **depressed quartic:**

$$\xi^4 + \left(3\sqrt{2} - \frac{11}{2}\right)\xi^2 + \left(4\sqrt{2}-6\right)\xi + \frac{33}{16} - \frac{5}{4}\sqrt{2} = 0 \tag{1.24}$$

And denoting the coefficients of this equation by

$$C \equiv 3\sqrt{2} - \frac{11}{2},\ D \equiv 4\sqrt{2} - 6,\ E \equiv \frac{33}{16} - \frac{5}{4}\sqrt{2} \tag{1.25}$$

the equation (1.24) can be written as,

$$\xi^4 + C\xi^2 + D\xi + E = 0 \tag{1.26}$$

Completing a square on the left side and transposing the other terms to the right side, we have,

$$\xi^4 + C\xi^2 + D\xi + E = 0 \Leftrightarrow \xi^4 + 2C\xi^2 + C^2 = C\xi^2 - D\xi + C^2 - E,$$

$$(\xi^2 + C)^2 = C\xi^2 - D\xi + C^2 - E \tag{1.27}$$

Let us transform (1.27) by introducing the variable $\zeta$ thus,

$$\left[(\xi^2+C)+\zeta\right]^2=\underbrace{C\xi^2-D\xi+C^2-E}_{(\xi^2+C)^2}+2\zeta(\xi^2+C)+\zeta^2 \tag{1.28}$$

This is the crucial trick, since for some $\zeta\in\mathbb{R}$ the right side of (1.28), which is quadratic in $\xi$, will be a *perfect square*, and the *necessary and sufficient condition* for this is that the *discriminant* of

$$(C+2\zeta)\xi^2-D\xi+C^2-E+2C\zeta+\zeta^2=0 \tag{1.29}$$

is null,

$$\Delta=0\Leftrightarrow(-D)^2-4(C+2\zeta)\left[C^2-E+2C\zeta+\zeta^2\right]=0 \tag{1.30}$$

cubic equation in $\zeta$:

$$\zeta^3+\frac{5C}{2}\zeta^2+(2C^2-E)\zeta+\left\{\frac{C(C^2-E)}{2}-\frac{D^2}{8}\right\}=0 \tag{1.31}$$

Denoting

$$c\equiv\frac{5C}{2},\ d\equiv2C^2-E,\ e\equiv\frac{C(C^2-E)}{2}-\frac{D^2}{8} \tag{1.32}$$

i.e.,

$$c\equiv\frac{5(6\sqrt{2}-11)}{4},\ d\equiv\frac{1511}{16}-\frac{259\sqrt{2}}{4},\ e\equiv\frac{5203\sqrt{2}}{32}-\frac{14769}{64} \tag{1.33}$$

we can write (1.31) in a more compact form

$$\zeta^3+c\zeta^2+d\zeta+e=0 \tag{1.34}$$

and applying again a *Tschirnhausen transformation*,

$$\zeta\equiv z-\frac{c}{3} \tag{1.35}$$

is eliminated the quadratic term and we get the **depressed cubic:**

$$z^3+\left(4\sqrt{2}-\frac{73}{12}\right)z+\left\{\frac{22\sqrt{2}}{3}-\frac{1133}{108}\right\}=0 \tag{1.36}$$

Denoting its coefficients by

$$\gamma\equiv4\sqrt{2}-\frac{73}{12},\ \delta\equiv\frac{22\sqrt{2}}{3}-\frac{1133}{108} \tag{1.37}$$

we can write (1.36) simply as

$$z^3+\gamma z+\delta=0 \tag{1.38}$$

To solve (1.38), let us make the change,

$$z \equiv u + v \tag{1.39}$$

Then, $z^3 + \gamma z + \delta = 0 \Leftrightarrow (u+v)^3 + \gamma(u+v) + \delta = 0 \Leftrightarrow u^3 + v^3 + 3uv(u+v) + \gamma(u+v) + \delta = 0$, i.e.,

$$u^3 + v^3 + (\gamma + 3uv)(u+v) + \delta = 0 \tag{1.40}$$

If we choose $u$, $v$ so that $\gamma + 3uv = 0$, we get two equations:

$$u^3 + v^3 = -\delta \tag{1.41}$$

$$u^3 v^3 = -\left(\frac{\gamma}{3}\right)^3 \tag{1.42}$$

Equations that give us the sum and product, respectively, of the cubes of $u$ and $v$, which allows us to construct a quadratic equation whose roots are these cubes:

$$t^2 + \delta t - \left(\frac{\gamma}{3}\right)^3 = 0 \tag{1.43}$$

The solutions of (1.43) are,

$$t = -\frac{\delta}{2} \pm \sqrt{\left(\frac{\delta}{2}\right)^2 + \left(\frac{\gamma}{3}\right)^3} \tag{1.44}$$

With no loss of generality, we can write:

$$u^3 = -\frac{\delta}{2} + \sqrt{\left(\frac{\delta}{2}\right)^2 + \left(\frac{\gamma}{3}\right)^3}, \quad v^3 = -\frac{\delta}{2} - \sqrt{\left(\frac{\delta}{2}\right)^2 + \left(\frac{\gamma}{3}\right)^3} \tag{1.45}$$

Thus, by (1.39), a solution of (1.38) is given by (*Cardano - Tartaglia formula*):

$$z = \sqrt[3]{-\frac{\delta}{2} + \sqrt{\left(\frac{\delta}{2}\right)^2 + \left(\frac{\gamma}{3}\right)^3}} + \sqrt[3]{-\frac{\delta}{2} - \sqrt{\left(\frac{\delta}{2}\right)^2 + \left(\frac{\gamma}{3}\right)^3}} \in \mathbb{R} \tag{1.46}$$

Note: this root is real for being $\gamma \geq 0$.

According to (1.37), the value of the discriminant of (1.38) is,

$$\Delta \equiv \left(\frac{\delta}{2}\right)^2 + \left(\frac{\gamma}{3}\right)^3 = \frac{2639}{108} - \frac{311\sqrt{2}}{18} \tag{1.47}$$

Let us write (1.46) in a more compact form,

$$z = \sqrt[3]{-\frac{\delta}{2} + \sqrt{\Delta}} + \sqrt[3]{-\frac{\delta}{2} - \sqrt{\Delta}} \tag{1.48}$$

Undoing now the variable change (1.35) and the value of *c* given by (1,33),

$$\zeta=\frac{5(11-6\sqrt{2})}{12}+\sqrt[3]{-\frac{\delta}{2}+\sqrt{\Delta}}+\sqrt[3]{-\frac{\delta}{2}-\sqrt{\Delta}} \qquad (1.49)$$

For the value of $\zeta$ given by (1.49) the quadratic equation (1.29) has a double root (since its discriminant vanishes). Therefore, (1.28) can be expressed as,

$$\left[(\xi^2+C)+\zeta\right]^2=\left\{\sqrt{C+2\zeta}\left[\xi-\frac{D}{2(C+2\zeta)}\right]\right\}^2 \qquad (1.50)$$

Extracting the square root of both sides of (1.50) and simplifying the right side,

$$\xi^2+C+\zeta=\pm\left(\xi\sqrt{C+2\zeta}-\frac{D}{2\sqrt{C+2\zeta}}\right) \qquad (1.51)$$

From (1.51) we get two quadratic equations in $\xi$,

$$\xi^2-\left\{\sqrt{C+2\zeta}\right\}\xi+\frac{D}{2\sqrt{C+2\zeta}}+C+\zeta=0 \qquad (1.52)$$

$$\xi^2+\left\{\sqrt{C+2\zeta}\right\}\xi-\frac{D}{2\sqrt{C+2\zeta}}+C+\zeta=0 \qquad (1.53)$$

and each of these ones provides a pair of solutions of the depressed quartic (1.26):

$$\xi=\frac{1}{2}\left[+\sqrt{C+2\zeta}\pm\sqrt{-\frac{2D}{\sqrt{C+2\zeta}}-3C-2\zeta}\right] \qquad (1.54)$$

$$\xi=\frac{1}{2}\left[-\sqrt{C+2\zeta}\pm\sqrt{+\frac{2D}{\sqrt{C+2\zeta}}-3C-2\zeta}\right] \qquad (1.55)$$

And finally, undoing the change of variable (1.23) and the value of *D*, given by (1,25),

$$\varepsilon=\frac{1}{2}\left[1-\sqrt{2}+\sqrt{C+2\zeta}\pm\sqrt{+\frac{4(3-\sqrt{2})}{\sqrt{C+2\zeta}}-3C-2\zeta}\right]\Rightarrow\begin{cases}\boxed{\varepsilon\approx 0.9476620415}\\ \varepsilon\approx 0.1762492425\end{cases} \qquad (1.56)$$

$$\varepsilon=\frac{1}{2}\left[1-\sqrt{2}-\sqrt{C+2\zeta}\pm\sqrt{-\frac{4(3-\sqrt{2})}{\sqrt{C+2\zeta}}-3C-2\zeta}\right]\Rightarrow\varepsilon\approx -0.9761692044\pm i\,0.2726245316 \qquad (1.57)$$

where the values of *C*, $\delta$, $\Delta$ and $\zeta$ are given by (1.25), (1.37), (1.47) and (1.49), respectively.

Obviously, in the context of the problem the last two roots has no sense, and the second, by the discussion made above, it is not valid. So, the only acceptable solution is the first, i.e., $r_1=r_2=r_3 \Leftrightarrow \varepsilon=0.94766204...$

**Relationship between *b* and *r*.**

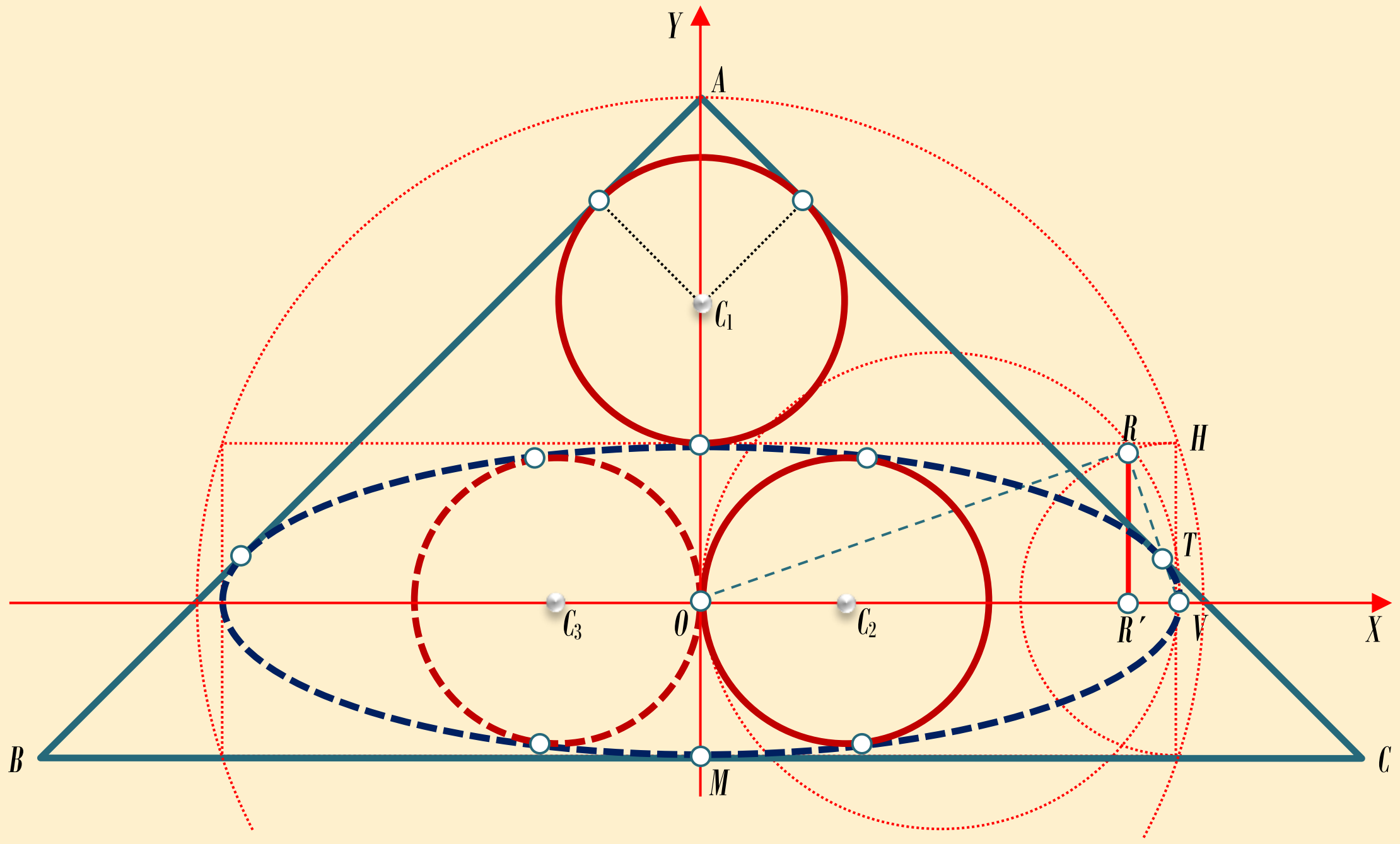

Straight from the above picture, we get the height $\overline{AM}$ of the triangle $\triangle CAB$:

$$\overline{AM} = 2\beta + r + r\sqrt{2} = 2\beta + r(1+\sqrt{2}) \quad (1.58)$$

Since $\triangle$AMC is an isosceles right triangle, the *Pythagorean Theorem* gives:

$$b = \overline{AC} = \overline{AM}\sqrt{2} = \left[2\beta + r(1+\sqrt{2})\right]\sqrt{2} = \left[2\frac{r}{\varepsilon} + r(1+\sqrt{2})\right]\sqrt{2} = r\left[2+\left(1+\frac{2}{\varepsilon}\right)\sqrt{2}\right],$$

i.e.,

$$\frac{b}{r} = 2 + \left(1+\frac{2}{\varepsilon}\right)\sqrt{2} \quad (1.59)$$

Therefore, an **explicit and algebraically exact relation between *b* and *r*** is

$$\frac{b}{r} = 2+\sqrt{2}+\frac{4\sqrt{2}}{1-\sqrt{2}+\sqrt{C+2\zeta}+\sqrt{\frac{4(3-2\sqrt{2})}{\sqrt{C+2\zeta}}-3C-2\zeta}} \quad (1.60)$$

where,

$$C \equiv 3\sqrt{2}-\frac{11}{2} \qquad \left.\begin{array}{l} \delta \equiv \dfrac{22\sqrt{2}}{3}-\dfrac{1133}{108} \\ \Delta \equiv \dfrac{2639}{108}-\dfrac{311}{9\sqrt{2}} \end{array}\right\} \qquad \zeta \equiv 5\left(\frac{11}{12}-\frac{1}{\sqrt{2}}\right)-\sum_{k=1}^{2}\sqrt[3]{\frac{\delta}{2}+(-1)^k\sqrt{\Delta}} \quad (1.61)$$

**Graphic representation of *r*.**

For a geometrically exactly representation of the circle inscribed in the right semi-ellipse, its radius can be determined as shown in the above drawing, where *R* is the intersection point of the circle with center at the midpoint of *OV* and diameter $\alpha$ and the circle with center at *V* and diameter $\beta$. Then,

$$\overline{RR'} \perp \overline{OV} \wedge \boxed{\overline{RR'} = r} \quad (1.62)$$

Indeed, $\triangle VRO$ is rectangle in *R*, and verified:

$$\left.\begin{array}{l} \overline{VO} = \alpha \\ \overline{VR} = \beta \\ \overline{RO} = \gamma \\ \triangle RR'O \sim \triangle VRO \end{array}\right\} \Rightarrow \frac{\overline{RR'}}{\overline{RO}} = \frac{\overline{VR}}{\overline{VO}} \Leftrightarrow \overline{RR'} = \overline{VR} \cdot \frac{\overline{RO}}{\overline{VO}} = \beta \frac{\gamma}{\alpha} = \beta\varepsilon = r \Leftrightarrow \boxed{\overline{RR'} = r} \quad (1.63)$$

Approximate value of the ratio $b/r$ with 3000 decimal digits:

$\frac{b}{r} \approx$ 6.3988504908313964106495157287278877674472204061565403354592771240691523730333207501758409383551522944
1023878174854651640468456996687068973239665149295906421075063249576063491593420540113369980768041554
2261184792106475817952240280238929292177543863534381751892191381811825172495071286753477187798972172
5157840494054209712156351268956764262678951334582589824139061257397271562674430007374644245088738121
5752937661619158810371550353830298896131654523806268076416303046466932699733493430716676084742789500
9321302367533275855874744625607455712042044600666498673003616785995329545391667895140857356972483175
5294249333962306618211976144376198184981958938277166998343743252036900095310115948025152210958909222
1513138042777850780756984606806710439625975680568147193395352556329568874837557216897273332382622548
0371780358763405280416802464854623121859352005127742007353166957344482980485700299651594172045582147
4873077788539077870006266513811650783029431840886716647639507358724606203091986829847162154480347058
9476760636794790961304363248728708893020273144713559708937065439810106019844467917147592222129386912
9613909871405662408940977216630807362222889432027339408862469931364161263733989411449518516228389543
3363979559022865920408475892931466073580359708382974248264028043891056391748424989408908598489895603
0043795034771317969424798623283799287907377834104974925525639573178304695558081124679143232290092016
4803997086156774871776646989224968250980134695323814202461067377099805287080637574741527617131082095
9743474016032491620666824325876162181374142222082275229416481573223019410147973230277950849712527959
7830384749181436162719413011813276117818766252675997809270813269283682912210223918787358216454740996
9528125748485094990366095890173191838947965469178100470935912783159420416570581077832288339098224207
3110575135847441526791042642455519302753735742708324551179866442909818824377480301585491913826110362
6079733946048440586666418704535999057312641290278938157958800034280640875230344487213731901329897046
5271245678278548364275127059546451323111360854817324731492716079419845748750894848355446064886710035
8238380604958221316342681631528164522697754746387028380248988448786159321388355504041402027009977829
1038149436337925696682679957368103700525421557030227133450494195164523937591121018859466812875691744
7673380426617813425144134488956013084158403250093059526778746910900957074848657666033251666524207326
6458808747162628561407010698819283256778346445819210020468153735707096776364173437639413662033632733
7804972836346149954942709839019424650587396671128037202704605746007755860203404091248831675040291503
5351193625155841623360353902411603638890965785836924747637039608669547026465316842958869552478543757
4608002133565052068045895206261961174996235896471178918631781434403967650673481708023304407049813951
7405160481488849112313246042656944956199015815721614276268625559567565298309106573451444981884817193
0485350645589954067006985105190995417008627930370132541628633310917934959265391196957533438702115954

## §2. Asymmetric case: $b \neq c$.

### 2.1. Existence and uniqueness of solution.

By observation of the below picture can be seen straightforward that *the solution is unique* in each quadrant, since the radius of the circles inscribed in the ellipse decreases as the eccentricity of the ellipse increases, according to (1.15) and (1.17), while the radius of the outside circle grows. Therefore, can only be a value of the eccentricity for which the radius of the three circles are equal. Hence, there is a bijection $\varepsilon \leftrightarrow r$.

On the other hand, the eccentricity increases when decreasing the angle $\varphi$ and, consequently, the radius $r$ increases as $\varphi$ increases. There is, therefore, also a bijection $\varphi \leftrightarrow r$, that applies $(90°, 0) \leftrightarrow (r_{max}, 0)$.

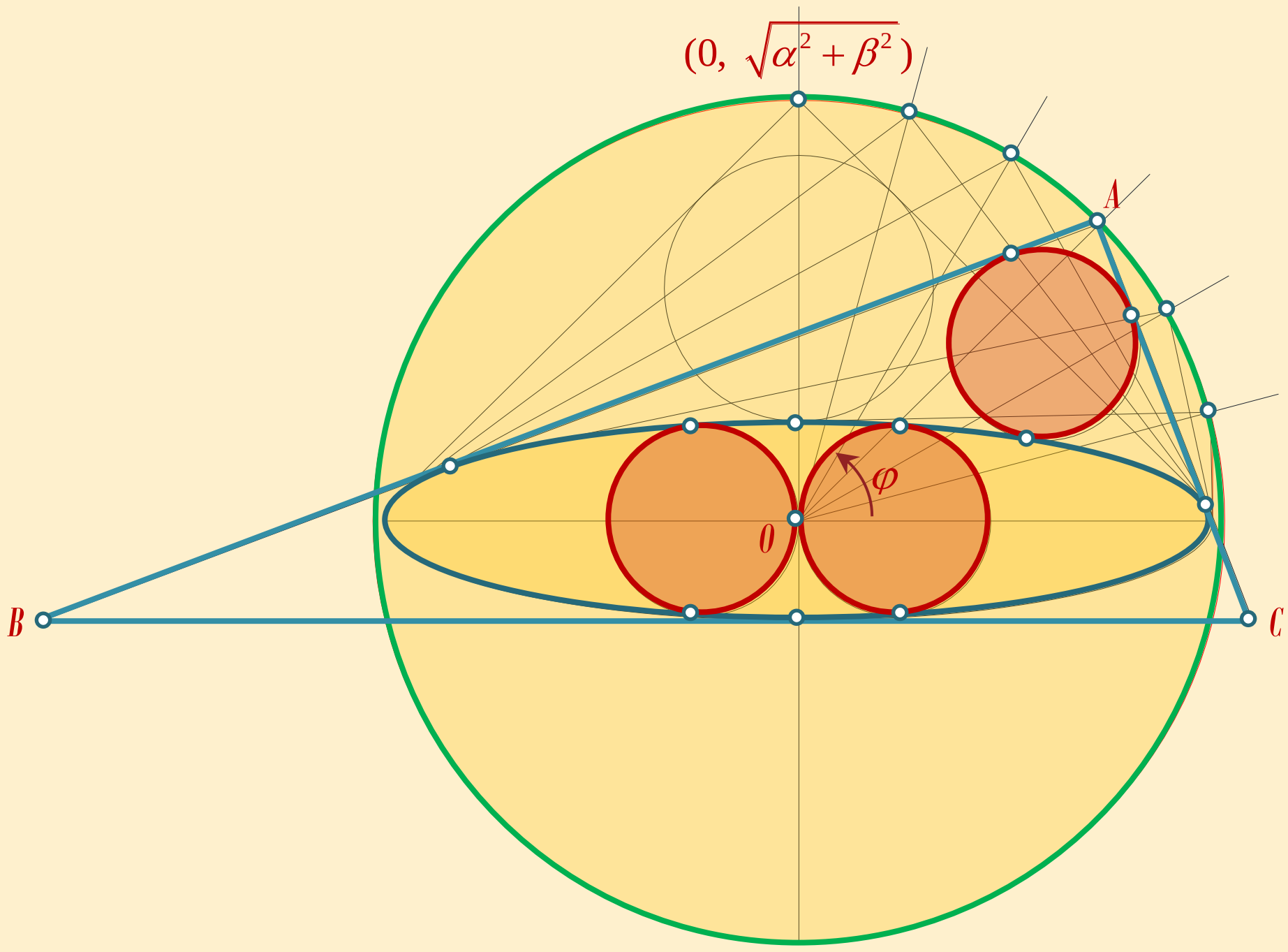

The radius for the three circles inscribed reaches a *maximum* when the value of eccentricity is the one obtained for the particular case studied before ($\varphi = 90° \Rightarrow \triangle ABC$ is a *right-angled isosceles triangle*).

But now the $\triangle ABC$ has no axial symmetry, so it is expected a significant complication in the calculations. The outline of several approaches (shown in appendixes B, C, D and F) from which outcrops an intractable complexity, allows glimpsing the impossibility of obtaining an explicit algebraic solution.

Despite that the extremely complexity of the relationship gotten for the symmetric case, doesn´t represents a mathematical proof of the impossibility of getting relationships for the asymmetric case, indeed it is a proof for the intuition and common sense: the asymmetry increases the complication, both in the outlines as in the difficult of the equations. Nevertheless, the aforementioned impossibility will be rigorously translated to an algebraic condition: the solubility of a sextic, (2.33).

## 2.2. Outlining the path to the solution: What we need for defining algebraically the geometric conditions?.

First of all, we must see the problem in a *bird view*, looking for the necessary equations for translate the geometric picture into an algebraic description.

In the study of the symmetric case we had yet determined the equation (1.13), that describe the *condition of tangency* imposed to the circumferences $\mathcal{C}_2$ and $\mathcal{C}_3$, and this relationship is also valid for the asymmetric case, since it doesn´t depend on $\triangle ABC$.

Other condition is that the ellipse $\mathcal{E}$ has to be tangent to the three sides of $\triangle ABC$. The *condition of tangency* among $\mathcal{E}$ and the hypotenuse $a$ is trivial for being parallel to $a$ the major axis of $\mathcal{E}$. The condition of tangency among $\mathcal{E}$ and the other two sides, $b$ and $c$, will be determine cutting $\mathcal{E}$ with the lines that contains $b$ and $c$, and imposing that the intersection verifies in a single point (the intersection of a line and an ellipse is determined by a quadratic equation, and the condition of tangency implies that its discriminant is null).

The *condition of tangency* between $\mathcal{C}_1$ and $b$ and $c$ is given simply imposing that the center $C_1$ of $\mathcal{C}_1$ belongs to the bisector line of the right-angle $A$.

The *condition of tangency* between $\mathcal{C}_1$ and $\mathcal{E}$ can be given imposing that, at an undefined point $T$ of $\mathcal{C}_1$ and $\mathcal{E}$, both curves share the tangent line, and this will occur iff the derivatives of $\mathcal{C}_1$ and $E$ take the same value at $T$ (i.e., tangents equal slope).

The equation of $\mathcal{E}$ involves two parameters, $\alpha$ and $\beta$, and the relationship (1.13) a third parameter, $\varepsilon$. Therefore, we need, at least, four independent equations for removing these three parameters. But by the definition of *eccentricity of the ellipse* we have an additional relationship between the three parameters involved, and by mean of this one and (1.13) the remaining equation can be easily reduced to a mono - parametric form, for instance, to an $\varepsilon$ - parametric for. It is necessary to find a third relationship, an $\varepsilon$ - parametric equation independent of the coordinates of $T$, for expressing $\varepsilon$ in the set of the final variables, $\mathcal{F} = \{a, b, c, r\}$. Let´s denote this hypothetic relationship as $\varepsilon_F$.

Substituting now $\varepsilon_F$ into the equation derived of the *condition of tangency* between $\mathcal{C}_1$ and $\mathcal{E}$, we´ll find an expression for one coordinate of T, (for instance, the abscise $x_T$) in terms of the FVs, and using the equation of $\mathcal{C}_1$ we could find the other coordinate of T (for instance, the ordinate $y_T$) in terms of the FVs.

Finally, equaling to $r$ the distance between $C_1$ and $T$, we would find the searched relationship among the variables of $\mathcal{F}$.

Nevertheless, we´ll see that the hardest difficulty is to find an explicit relationship $\varepsilon_F$, that would be the key for the door to the solution, but, as it will be shown later, $\varepsilon_F$ is, in fact, the key for proving that it is impossible to find a relation in the set $\mathcal{F} = \{a, b, c, r\}$, neither explicit nor implicit, since the equation relating $r$ to $\varepsilon$ is of 6th grade in $\varepsilon$, without any especial relation of symmetry in its coefficients, so unsolvable by radicals.

However, taking as final set of variables $\mathcal{F}' = \{a, b, c, \varepsilon\}$, it is indeed possible establish a relationship, because the equation relating $r$ to $\varepsilon$ is of 4th grade in $r$, therefore solvable algebraically. Getting a relationship in $\mathcal{F}' = \{a, b, c, \varepsilon\}$ allows to get another one in $\mathcal{F} = \{a, b, c, r\}$ for a given right-angled triangle $\triangle ABC$, by means of numeric methods.

## 2.3. Switching the reference system.

The previous reference system provides the simplest expression for the equation of the ellipse (reduced to its canonical form), but instead, the equations for the remaining elements result too complex (in appendix B we show its expressions). Let us take now the Cartesian reference system showed in the figure below, obtained from the previous one by a rotation and a translation of axes: rotation of 180º (keeping the parallelism of the coordinate's axes with respect to the axis of the ellipse) and translation of the origin to the vertex *A*.

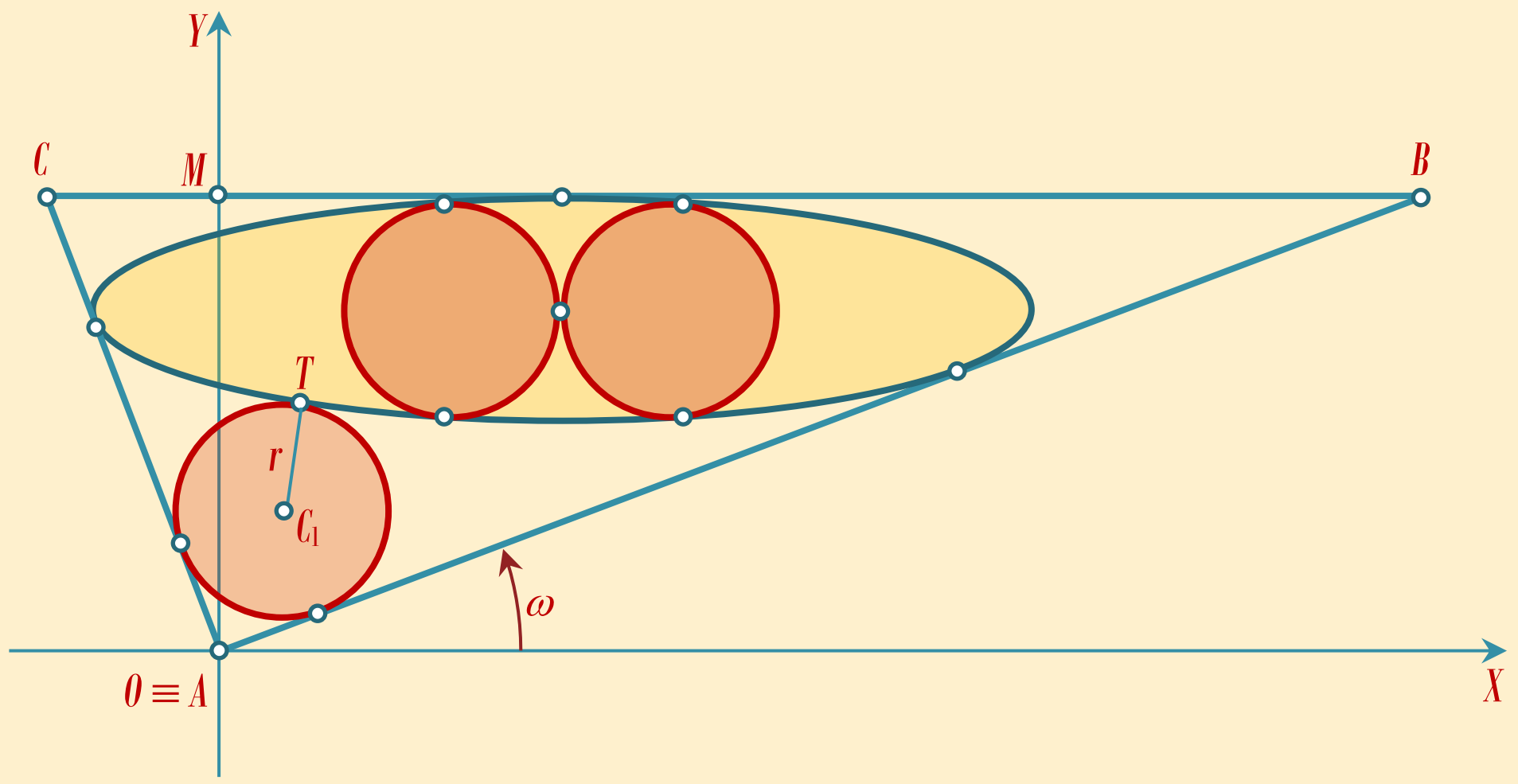

## 2.3. Equations in the new reference system.

To relate the new coordinates with the old, just switch $x$ by $-x + h$, and the primitives coordinates of $C_1$ with the new ones. The slope *m* is invariant under the transformation performed (180 degree turn and translation). In this new scenery, the equation of the bisector of A, and the equations for the coordinates of $C_1$, are significantly simplify.

The slope of the bisector of *A* is,

$$\tan\left(\omega+\frac{\pi}{4}\right)=\frac{\tan\omega+\tan\frac{\pi}{4}}{1-\tan\omega\tan\frac{\pi}{4}}=\frac{1+\tan\omega}{1-\tan\omega},$$

where $\tan\omega$ is the slope of the side *AB*. Setting,

$$\boxed{m\equiv\tan\omega} \qquad (2.01)$$

the equation of the bisector of *A* is as simply as

$$\boxed{y=\frac{1+m}{1-m}x} \qquad (2.02)$$

The new coordinates of $C_1$ can be obtained as intersection of a circumference with center at the origin and radius equal to the diagonal of a square of side *r*. Assuming $m > 0$,

$$\left.\begin{array}{l} y=\dfrac{1+m}{1-m}x \\ x^2+y^2=(r\sqrt{2})^2 \end{array}\right\}\Rightarrow\left[1+\left(\frac{1+m}{1-m}\right)^2\right]x^2=2r^2$$

$$\boxed{x_1=\frac{|1-m|\,r}{\sqrt{1+m^2}}}\quad\boxed{y_1=\frac{(1+m)r}{\sqrt{1+m^2}}} \qquad (2.03)$$

and expressing $m$ as a function of the sides of $\triangle ABC$,

$$x_1 = \frac{|c-b|r}{a} \qquad y_1 = \frac{(c+b)r}{a} \tag{2.04}$$

## 2.4. Centre of the ellipse: analytic and geometric determination.

The equation of the ellipse with center at $E(x_0, y_0)$ is,

$$\mathscr{E} \equiv \frac{(x-x_0)^2}{\alpha^2} + \frac{(y-y_0)^2}{\beta^2} = 1 \tag{2.05}$$

and the equation of the line determined by the vertex of $A$ and $B$,

$$y = mx \tag{2.06}$$

The intersection of this line and the ellipse is given by,

$$\frac{(x-x_0)^2}{\alpha^2} + \frac{(mx-y_0)^2}{\beta^2} = 1 \Leftrightarrow (x-x_0)^2 + (mx-y_0)^2 - \alpha^2\beta^2 = 0,$$

or,

$$(m^2\alpha^2 + \beta^2)x^2 - 2(\beta^2 x_0 + m\alpha^2 y_0)x + \beta^2 x_0 + \alpha^2 y_0 - \alpha^2\beta^2 = 0 \tag{2.07}$$

but the line (2.06) will be tangent to the ellipse (2.05) iff the discriminant of the quadratic equation in $x$ (2.07) is null,

$$\Delta = \left[-2(\beta^2 x_0 + m\alpha^2 y_0)\right]^2 - 4(m^2\alpha^2 + \beta^2)(\beta^2 x_0 + \alpha^2 y_0 - \alpha^2\beta^2) = 0 \tag{2.08}$$

and simplifying,

$$m^2(x_0^2 - \alpha^2) - 2mx_0y_0 + y_0^2 - \beta^2 = 0 \tag{2.09}$$

Substituting in (2.09) $m$ by $-1/m$ we obtain a similar expression for the condition of tangency of the line determined by the vertices $A$ and $C$ with the ellipse (2.05):

$$x_0^2 - \alpha^2 + 2mx_0y_0 + m^2(y_0^2 - \beta^2) = 0 \tag{2.10}$$

Adding equations (2.09) and (2.10) side by side, and factoring, we have,

$$(m^2+1)(x_0^2 + y_0^2 - \alpha^2 - \beta^2) = 0 \tag{2.11}$$

that implies,

$$\alpha^2 + \beta^2 = x_0^2 + y_0^2 \tag{2.12}$$

This condition is also obtained considering that the vertex $A$ belongs to the ***director circumference*** (see Appendix B).

From (2.12),

$$\alpha^2 = x_0^2 + y_0^2 - \beta^2 \tag{2.13}$$

and substituting $\alpha$ in (2.09) or (2.10) by the right side of (2.13),

$$(y_0^2 - \beta^2)(1 - m^2) = 2 m x_0 y_0 \tag{2.14}$$

Into this equation hides an interesting geometric relationship, which comes to light when expressed thus,

$$\frac{2m}{1-m^2} = \frac{y_0^2 - \beta^2}{x_0 y_0} \tag{2.15}$$

since,

$$\frac{2m}{1-m^2} = \tan 2\omega \tag{2.16}$$

Using (2.16) and writing right now (2.14) in the shape,

$$\frac{x_0}{y_0 - \beta} = \frac{(y_0 + \beta)\cot 2\omega}{y_0} \tag{2.17}$$

follows a useful geometrical construction to determine the center of the ellipse in terms of $\beta$:

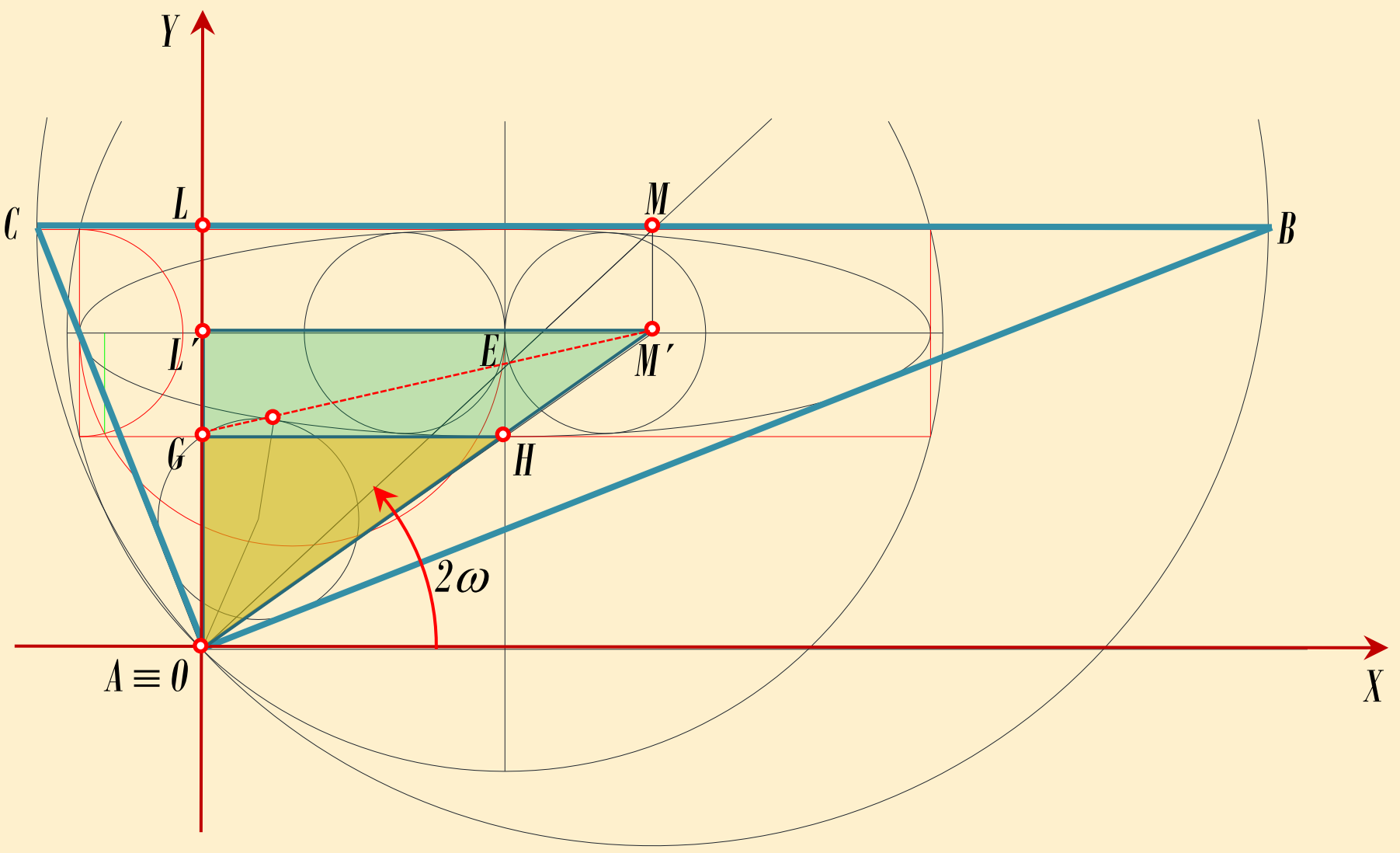

$$\frac{\overline{GH}}{\overline{OG}} = \frac{\overline{L'M'}}{\overline{OL'}} = \frac{\overline{LM}}{\overline{OL'}} \Leftrightarrow \frac{x_0}{y_0 - \beta} = \frac{(y_0 + \beta)\cot 2\omega}{y_0} \tag{2.18}$$

But $\triangle OBM'$ and $\triangle OCM'$ are isosceles because both triangles have two equals angles; hence, $M$ is the middle point of the hypotenuse of $\triangle ABC$, and the vertex $H$ of the ellipse is simply the intersection of $AM'$ with the parallel to $BC$ by $G$, where $\overline{GL} = 2\beta$ and $M'$ is the projection of $M$ over the major axis of the ellipse.

From (2.17), taking into account that $y_0+\beta$ is the height of $\triangle ABC$ relative to its hypotenuse, and can be expressed as,

$$y_0+\beta=\frac{bc}{a} \tag{2.19}$$

since $m$ is the slope of $c$ in $\triangle ABC$ *and* $m=\tan\hat{B}=b/c$.

Solving for $x_0$ in (2.17) and replacing $y_0+\beta$ by the value given by (2.19), since $\tan\,\omega=b/c$,

$$x_0=[(y_0+\beta)\cot 2\omega]\cdot\frac{y_0-\beta}{y_0}=\frac{bc}{a}\cdot\frac{c^2-b^2}{2bc}\left(1-\frac{\beta}{y_0}\right),$$

$$x_0=\frac{c^2-b^2}{2a}\left(1-\frac{\beta}{y_0}\right) \tag{2.20}$$

Note that, $\dfrac{c^2-b^2}{2a}=\dfrac{-\frac{b^2}{a}+\frac{c^2}{a}}{2}=\dfrac{x_C+x_B}{2}=x_M$, that is, the *abscise of the middle point* of $\overline{BC}$.

Let´s denote,

$$k\equiv\frac{c^2-b^2}{2a} \tag{2.21}$$

$$h\equiv\frac{bc}{a} \tag{2.22}$$

**Theorem.** The *locus of the centers of the ellipses* with minor semi-axis of length $\beta$ *inscribed in* a right-angled triangle $\triangle ABC$, with major-axis parallel to the hypotenuse, in a Cartesian reference system with origin at the vertex $A$ of the right-angle $\beta$ and abscissa axis below the hypotenuse and parallel to it, belong to the ***hyperbole***,

$$y=\frac{\beta}{1-\frac{x}{k}} \tag{2.23}$$

Removing $y_0$ in (2.20) by means of (2.19), we obtain the **$\beta$-parametric equations of the center of the ellipse:**

$$\begin{cases} x_0=k\left[1-\left(\dfrac{h}{\beta}-1\right)^{-1}\right] \\ y_0=h-\beta \end{cases} \tag{2.23}$$

and by eliminating $\beta$ through (1.13) we get the **$\varepsilon$-parametric equations of the center of the ellipse:**

$$\begin{cases} x_0=k\left(1-\dfrac{r}{h\varepsilon-r}\right) \\ y_0=h-\dfrac{r}{\varepsilon} \end{cases} \tag{2.24}$$

Notice the simplified relationship,

$$h^2+k^2=\left(\frac{a}{2}\right)^2 \tag{2.25}$$

## 2.5. Coordinates of the tangent point $T$ between $\mathcal{C}_1$ and $\mathcal{E}$.

For determining the point of tangency between the $\mathcal{C}_1$ and $\mathcal{E}$, as shown below, it is possible to raise three independent equations. At first glance, defining in three different ways the same object (unknown $x_T$) might seem unnecessary work, and therefore superfluous. However, is quite the opposite, since the implementation of the TOP (which displays here some of its potential) reduces the three quartic to a simple quadratic equation which preserves the root of interest (the abscissa of $T$) shared by the three quartics. The TOP acts synergistically on the equations, so that the joint action of the equations is much more effective than just the sum of their individual actions. Below are defined the three quartics.

### 2.5.1. Quartic for $x_T$ derived from the intersection of $\mathcal{C}_1$ and $\mathcal{E}$.

$$\left.\begin{array}{l}\dfrac{(x-x_0)^2}{\alpha^2}+\dfrac{(y-y_0)^2}{\beta^2}=1\Leftrightarrow y=y_0\pm\beta\sqrt{1-\dfrac{(x-x_0)^2}{\alpha^2}}\\ (x-x_1)^2+(y-y_1)^2=r^2\Leftrightarrow y=y_1\pm\sqrt{r^2-(x-x_1)^2}\end{array}\right\}\Rightarrow y_0\pm\beta\sqrt{1-\dfrac{(x-x_0)^2}{\alpha^2}}=y_1\pm\sqrt{r^2-(x-x_1)^2}.$$

Taking into account that the contact take place between the *lower semi-ellipse* and the *upper semi-circumference*, the signs of the radicals must be $-$ in the equation of the ellipse and $+$ in the equation of the circumference:

$$y_0-\beta\sqrt{1-\frac{(x-x_0)^2}{\alpha^2}}=y_1+\sqrt{r^2-(x-x_1)^2} \qquad (2.26)$$

Streamlining, developing and collecting terms (used the DERIVE computer algebra system), there arises a quartic equation, that we denote by $u(x)=0$, showing the following dissuasive-looking: (2.27)

$$u(x)=a^2 ε^6 (h ε-r)^4 x^4 +4 a ε^4 (r-h ε)^3 (a k (ε^2 -1) (h ε-2 r)+r (b-c) (r-h ε)) x^3 -2 (a^2 (h^4 ε^4 (ε^2 -2)+4 h^3 r ε^3 (2-ε^2)+h^2 ε^2 (r^2 (ε^4 +5 ε^2 -12)-k^2 ε^2$$
$$(ε^2 -1) (3 ε^2 -2))+2 h r ε (2 k^2 ε^2 (ε^2 -1) (3 ε^2 -2)-r^2 (ε^4 +ε^2 -4))-r^2 (4 k^2 ε^2 (ε^2 -1) (3 ε^2 -2)+r^2 (2-ε^4)))+2 a r ε (r-h ε) (b (h^2 ε^2 (ε^2 -2)-2 h ε (k$$
$$ε (ε^2 -1)+r (ε^2 -2))+r (4 k ε (ε^2 -1)+r (ε^2 -2)))+c (h^2 ε^2 (ε^2 -2)+2 h ε (k ε (ε^2 -1)+r (2-ε^2))-r (4 k ε (ε^2 -1)+r (2-ε^2)))-4 r^2 ε^2 (b^2 -b c ε^2 +c^2) (h$$
$$ε-r)^2 ) (h ε-r)^2 x^2 +4 (a^2 k ε (ε^2 -1) (h ε-2 r) (h^4 ε^3 -4 h^3 r ε^2 +h^2 ε (r^2 (ε^2 +5)-k^2 ε^2 (ε^2 -1))+2 h r (2 k^2 ε^2 (ε^2 -1)-r^2 (ε^2 +1))-r^2 ε (4 k^2$$
$$(ε^2 -1)-r^2 ))+a r (b (h^4 ε^4 -2 h^3 ε^3 (k ε (ε^2 -1)+2 r)+h^2 ε^2 (k^2 ε^2 (ε^2 -1)+8 k r ε (ε^2 -1)+r^2 (ε^2 +7))-2 h r ε (2 k^2 ε^2 (ε^2 -1)+5 k r ε (ε^2 -1)+r^2$$
$$(ε^2 +3))+r^2 (4 k^2 ε^2 (ε^2 -1)+4 k r ε (ε^2 -1)+r^2 (ε^2 +2)))-c (h^4 ε^4 +2 h^3 ε^3 (k ε (ε^2 -1)-2 r)+h^2 ε^2 (k^2 ε^2 (ε^2 -1)+8 k r ε (1-ε^2 )+r^2 (ε^2 +7))-2 h r ε (2$$
$$k^2 ε^2 (ε^2 -1)+5 k r ε (1-ε^2 )+r^2 (ε^2 +3))+r^2 (4 k^2 ε^2 (ε^2 -1)+4 k r ε (1-ε^2 )+r^2 (ε^2 +2)))) (h ε-r)-2 r^2 ε (b^2 (h ε-r)^2 +2 b c k ε (1-ε^2 ) (h ε-2 r)-c^2 (h$$
$$ε-r)^2 ) (h ε-r)^2 ) (h ε-r) x+a^2 (h^8 ε^6 -8 h^7 r ε^5 -2 h^6 ε^4 (k^2 ε^2 (ε^2 -1)+r^2 (ε^2 -13))+4 h^5 r ε^3 (4 k^2 ε^2 (ε^2 -1)+r^2 (3 ε^2 -11))+h^4 ε^2 (k^4 ε^4 (ε^2 -1)^2 -2 k^2$$
$$r^2 ε^2 (ε^4 +24 ε^2 -25)+r^4 (5 ε^4 -32 ε^2 +41))-4 h^3 r ε (2 k^4 ε^4 (ε^2 -1)^2 -k^2 r^2 ε^2 (3 ε^4 +16 ε^2 -19)+r^4 (5 ε^4 -12 ε^2 +5))+2 h^2 r^2 (12 k^4 ε^4 (ε^2 -1)^2 -k^2 r^2 ε^2$$
$$(13 ε^4 +15 ε^2 -28)+r^4 (15 ε^4 -21 ε^2 +2))-4 h r^3 ε (8 k^4 ε^2 (ε^2 -1)^2 -2 k^2 r^2 (3 ε^4 -ε^2 -2)+5 r^4 (ε^2 -1))+r^4 (16 k^4 ε^2 (ε^2 -1)^2 +8 k^2 r^2 ε^2 (1-ε^2 )+r^4 (5$$
$$ε^2 -4)))+4 a r (b+c) (r-h ε)^3 (h^4 ε^3 -4 h^3 r ε^2 +h^2 ε (r^2 (ε^2 +5)-k^2 ε^2 (ε^2 -1))+2 h r (2 k^2 ε^2 (ε^2 -1)-r^2 (ε^2 +1))-r^2 ε (4 k^2 (ε^2 -1)-r^2 ))+4 r^2 (b^2 (2 h^2$$
$$ε^2 -4 h r ε+r^2 (2-ε^2 )) (h ε-r)^2 +2 b c ε (h^4 ε^3 -4 h^3 r ε^2 +h^2 ε (5 r^2 -k^2 ε^2 (ε^2 -1))+2 h r (2 k^2 ε^2 (ε^2 -1)-r^2 )+4 k^2 r^2 ε (1-ε^2 ))+c^2 (2 h^2 ε^2 -4 h r$$
$$ε+r^2 (2-ε^2 )) (h ε-r)^2 ) (h ε-r)^2 =0$$

### 2.5.2. Quartic for $x_T$ which is obtained by equating the derivatives of $\mathcal{C}_1$ and $\mathcal{E}$ in its explicit form.

**Derivatives** of the upper canonic *semi-ellipse* $\mathcal{E}$ and the lower *semi-circumference* $\mathcal{C}_1$:

$$\frac{(x-x_0)^2}{\alpha^2}+\frac{(y-y_0)^2}{\beta^2}=1\Leftrightarrow y=y_0\pm\beta\sqrt{1-\frac{(x-x_0)^2}{\alpha^2}}\Rightarrow y'=\mp\beta\frac{\dfrac{x-x_0}{\alpha^2}}{\sqrt{1-\dfrac{(x-x_0)^2}{\alpha^2}}},$$

$$y'=-\frac{\beta}{\alpha}\frac{1}{\sqrt{\left(\dfrac{\alpha}{x_0-x}\right)^2-1}} \qquad (2.28)$$

$$(x-x_1)^2+(y-y_1)^2=r^2 \Leftrightarrow y=y_1\pm\sqrt{r^2-(x-x_1)^2} \Rightarrow y'=\mp\frac{x-x_1}{\sqrt{r^2-(x-x_1)^2}},$$

$$y'=-\frac{1}{\sqrt{\left(\frac{r}{x-x_1}\right)^2-1}} \qquad (2.29)$$

the derivatives of both curves at the tangent point *T* represent the slope of the common tangent line; hence, they have to be equals:

$$-\frac{\beta}{\alpha}\frac{1}{\sqrt{\left(\frac{\alpha}{x_0-x}\right)^2-1}}=-\frac{1}{\sqrt{\left(\frac{r}{x-x_1}\right)^2-1}} \Leftrightarrow \frac{\beta}{\alpha}\sqrt{\left(\frac{r}{x-x_1}\right)^2-1}=\sqrt{\left(\frac{\alpha}{x_0-x}\right)^2-1}$$

$$\left(\frac{\beta}{\alpha}\right)^2\left[\left(\frac{r}{x-x_1}\right)^2-1\right]=\left(\frac{\alpha}{x_0-x}\right)^2-1 \qquad (2.30)$$

where $\alpha^2$ can be expressed in terms of *r* and $\varepsilon$,

$$\alpha^2=\frac{r^2}{\varepsilon^2(1-\varepsilon^2)} \qquad (2.31)$$

From (2.30) and (2.31),

$$1-\left(\frac{\beta}{\alpha}\right)^2=\left(\frac{\alpha}{x_0-x}\right)^2-\left(\frac{\beta}{\alpha}\right)^2\left(\frac{r}{x-x_1}\right)^2 \Leftrightarrow \varepsilon^2=\frac{r^2}{\varepsilon^2(1-\varepsilon^2)(x_0-x)^2}-\frac{r^2(1-\varepsilon^2)}{(x-x_1)^2}$$

or,

$$\varepsilon^2(1-\varepsilon^2)(x_0-x)^2\left[\varepsilon^2(x-x_1)^2+r^2(1-\varepsilon^2)\right]-r^2(x-x_1)^2=0 \qquad (2.32)$$

Developing and collecting terms in x we obtain the following quartic denoted $v(x)=0$: (2.33)

$$v(x)=a^2 r^2 \varepsilon^4 (1-\varepsilon^2)(h\varepsilon-r)^2 x^4+2a\varepsilon^4(1-\varepsilon^2)(a k (h\varepsilon-2r)+r(b-c)(r-h\varepsilon))(h\varepsilon-r)x^3+(a^2(h^2\varepsilon^2(k^2\varepsilon^4(\varepsilon^2-1)-r^2(\varepsilon^6-2\varepsilon^4+\varepsilon^2-1))+2hr\varepsilon(r^2(\varepsilon^6-2\varepsilon^4+\varepsilon^2-1)-2k^2\varepsilon^4$$
$$(\varepsilon^2-1))+r^2(4k^2\varepsilon^4(\varepsilon^2-1)-r^2(\varepsilon^6-2\varepsilon^4+\varepsilon^2-1)))+4akr\varepsilon^4(1-\varepsilon^2)(b-c)(h\varepsilon-r)(h\varepsilon-2r)+r^2\varepsilon^4(b-c)^2(\varepsilon^2-1)(h\varepsilon-r)^2)x^2+2r(a^2kr\varepsilon^2(\varepsilon^2-1)^2(h\varepsilon-r)(h\varepsilon-2r)+a$$
$$(b-c)(h^2\varepsilon^2(k^2\varepsilon^4(\varepsilon^2-1)+r^2)-2hr\varepsilon(2k^2\varepsilon^4(\varepsilon^2-1)+r^2)+r^2(4k^2\varepsilon^4(\varepsilon^2-1)+r^2))+kr\varepsilon^4(1-\varepsilon^2)(b-c)^2(h\varepsilon-r)(h\varepsilon-2r))x-r^2(a^2k^2\varepsilon^2(\varepsilon^2-1)^2(h\varepsilon-2r)^2-(b-c)^2$$
$$(h^2\varepsilon^2(k^2\varepsilon^4(\varepsilon^2-1)+r^2)-2hr\varepsilon(2k^2\varepsilon^4(\varepsilon^2-1)+r^2)+r^2(4k^2\varepsilon^4(\varepsilon^2-1)+r^2)))=0$$

### 2.5.3. Quartic for $x_T$ got by equating the derivatives of $\mathcal{C}_1$ and $\mathcal{E}$ in its implicit form.

Deriving implicitly the equations of both curves,

$$\mathcal{E}\equiv\frac{(x-x_0)^2}{\alpha^2}+\frac{(y-y_0)^2}{\beta^2}=1\Rightarrow\frac{\not{2}(x-x_0)}{\alpha^2}+\frac{\not{2}(y-y_0)}{\beta^2}y'=0\Rightarrow y'=-\frac{\beta^2}{\alpha^2}\cdot\frac{x-x_0}{y-y_0} \qquad (2.34)$$

$$\mathcal{C}_1\equiv(x-x_1)^2+(y-y_1)^2=r^2\Rightarrow\not{2}(x-x_1)+\not{2}(y-y_1)y'=0\Rightarrow y'=-\frac{x-x_1}{y-y_1} \qquad (2.35)$$

But $\mathcal{C}_1$ and $\mathcal{E}$ share the tangent line at $T$, hence, their derivatives (slope of the tangent line) take the same value at $T$:

$$-\frac{\beta^2}{\alpha^2}\cdot\frac{x-x_0}{y-y_0}=-\frac{x-x_1}{y-y_1} \tag{2.36}$$

Parameterizing (2.36) in $\varepsilon$, eliminating denominators and expressing it in explicit form,

$$(1-\varepsilon^2)\cdot(x-x_0)(y-y_1)=(x-x_1)(y-y_0) \tag{2.37}$$

$$y=\frac{y_1(1-\varepsilon^2)(x-x_0)-y_0(x-x_1)}{(1-\varepsilon^2)\cdot(x-x_0)-(x-x_1)} \tag{2.38}$$

Replacing by means of (2.38) the variable $y$ in the equation of $\mathcal{C}_1$,

$$\mathcal{C}_1\equiv(x-x_1)^2+\left[\frac{y_1(1-\varepsilon^2)(x-x_0)-y_0(x-x_1)}{(1-\varepsilon^2)\cdot(x-x_0)-(x-x_1)}-y_1\right]^2=r^2 \tag{2.39}$$

$$\left[(x-x_1)^2-r^2\right]\left[\varepsilon^2x+x_0(1-\varepsilon^2)-x_1\right]^2+(y_0-y_1)^2(x-x_1)^2=0 \tag{2.40}$$

after a long chain of operations and simplifications, we get the following quartic that we denote by $w(x)=0$, (2.41)

$$w(x)=a^4 \varepsilon^6 (h \varepsilon-r)^2 x^4 +2 a^3 \varepsilon^4 (r-h \varepsilon) (a k (\varepsilon^2 -1) (h \varepsilon-2 r)+r (b-c) (\varepsilon^2 +1) (r-h \varepsilon)) x^3 +a^2 (a^2 (h^4 \varepsilon^4 -4 h^3 r \varepsilon^3 +h^2 \varepsilon^2 (k^2 \varepsilon^2 (\varepsilon^2 -1)^2 +r^2 (6-\varepsilon^6 ))-2 h$$
$$r \varepsilon (2 k^2 \varepsilon^2 (\varepsilon^2 -1)^2 +r^2 (2-\varepsilon^6 ))+r^2 (\varepsilon^2 -1) (4 k^2 \varepsilon^2 (\varepsilon^2 -1)-r^2 (\varepsilon^4 +\varepsilon^2 +1)))+2 a r \varepsilon (r-h \varepsilon) (b (h^2 \varepsilon^2 +h \varepsilon (k \varepsilon (2 \varepsilon^4 -\varepsilon^2 -1)-2 r)-r (2 k \varepsilon (2$$
$$\varepsilon^4 -\varepsilon^2 -1)-r))+c (h^2 \varepsilon^2 -h \varepsilon (k \varepsilon (2 \varepsilon^4 -\varepsilon^2 -1)+2 r)+r (2 k \varepsilon (2 \varepsilon^4 -\varepsilon^2 -1)+r)))+r^2 \varepsilon^2 (a^2 (\varepsilon^4 +4 \varepsilon^2 +2)-2 b c \varepsilon^2 (\varepsilon^2 +4)) (h \varepsilon-r)^2 ) x^2 +2 a r (a^3 k r \varepsilon^4$$
$$(\varepsilon^2 -1) (h \varepsilon-r) (h \varepsilon-2 r)+a^2 (b-c) (h^4 \varepsilon^4 -4 h^3 r \varepsilon^3 +h^2 \varepsilon^2 (k^2 \varepsilon^2 (\varepsilon^2 -1)^2 +r^2 (6-\varepsilon^4 ))-2 h r \varepsilon (2 k^2 \varepsilon^2 (\varepsilon^2 -1)^2 +r^2 (2-\varepsilon^4 ))+r^2 (\varepsilon^2 -1) (4 k^2 \varepsilon^2$$
$$(\varepsilon^2 -1)-r^2 (\varepsilon^2 +1)))+a r \varepsilon (b-c) (r-h \varepsilon) (b (2 h^2 \varepsilon^2 +h \varepsilon (k \varepsilon (\varepsilon^2 -1) (\varepsilon^2 +2)-4 r)-2 r (k \varepsilon (\varepsilon^2 -1) (\varepsilon^2 +2)-r))+c (2 h^2 \varepsilon^2 -h \varepsilon (k \varepsilon (\varepsilon^2 -1) (\varepsilon^2 +2)+4 r)+2$$
$$r (k \varepsilon (\varepsilon^2 -1) (\varepsilon^2 +2)+r)))+r^2 \varepsilon^2 (b-c) (b^2 (\varepsilon^2 +2)-2 b c \varepsilon^2 +c^2 (\varepsilon^2 +2)) (h \varepsilon-r)^2 ) x-r^2 (a^4 k^2 \varepsilon^2 (\varepsilon^2 -1)^2 (h \varepsilon-2 r)^2 +2 a^3 k r \varepsilon^2 (1-\varepsilon^2 ) (b-c) (h \varepsilon-r)$$
$$(h \varepsilon-2 r)-a^2 (b-c)^2 (h^2 \varepsilon^2 (h^2 \varepsilon^2 -4 h r \varepsilon+(k^2 \varepsilon^2 (\varepsilon^2 -1)^2 +r^2 (6-\varepsilon^2 )))-2 h r \varepsilon (2 k^2 \varepsilon^2 (\varepsilon^2 -1)^2 +r^2 (2-\varepsilon^2 ))+r^2 (\varepsilon^2 -1) (4 k^2 \varepsilon^2 (\varepsilon^2 -1)-r^2 ))+2 a r \varepsilon$$
$$(b-c)^2 (b (h^2 \varepsilon^2 +h \varepsilon (k \varepsilon (\varepsilon^2 -1)-2 r)-r (2 k \varepsilon (\varepsilon^2 -1)-r))+c (h^2 \varepsilon^2 -h \varepsilon (k \varepsilon (\varepsilon^2 -1)+2 r)+r (2 k \varepsilon (\varepsilon^2 -1)+r))) (h \varepsilon-r)-2 a^2 (b-c)^2 r^2 \varepsilon^2 (h \varepsilon-r)^2 )=0$$

or, more compact, (2.42)

$$w(x)=\varepsilon^4 x^4 -2 \varepsilon^2 x^3 (x_0 (\varepsilon^2 -1)+x_1 (\varepsilon^2 +1))-x^2 (r^2 \varepsilon^4 -x_0^2 (\varepsilon^2 -1)^2 -2 x_0 x_1 (2 \varepsilon^4 -\varepsilon^2 -1)-x_1^2 (\varepsilon^4 +4 \varepsilon^2 +1)-(y_0-y_1)^2 )+2 x (r^2 \varepsilon^2 (x_0 (\varepsilon^2 -1)+x_1)-x_1 (x_0^2$$
$$(\varepsilon^2 -1)^2 +x_0 x_1 (\varepsilon^2 -1) (\varepsilon^2 +2)+x_1^2 (\varepsilon^2 +1)+(y_0-y_1)^2 ))-r^2 (x_0^2 (\varepsilon^2 -1)^2 +2 x_0 x_1 (\varepsilon^2 -1)+x_1^2 )+x_1^2 (x_0^2 (\varepsilon^2 -1)^2 +2 x_0 x_1 (\varepsilon^2 -1)+x_1^2 +(y_0-y_1)^2 )=0$$

The quartic equations (2.27), (2.33) and (2.41) have to share the *real root* $x=x_T$ at the tangent point $T$. This property can be exploited applying the *Theorem of Overlapped Polynomials* (TOP), very useful for dealing with the enormous complexity of these equations. By applying the TOP three times, we get a quadratic from three quartics. First, applying the TOP to the quartics $u$ an $v$ it is obtained a cubic, denoted $\Delta_{uv}$, then, applying it to the quartics $v$ an $w$ it is obtained other cubic, denoted $\Delta_{vw}$. Finally, the application of the TOP to this pair of cubic equations it is obtained a quadratic equation, denoted $\Delta_{uvw}$. The algorithm can be schematized as follow: (2.43)

$$\left.\begin{array}{l} u,\ v\in\mathcal{O}_4(x-x_T)\Rightarrow\Delta_{uv}\equiv \mathrm{M}\{u\}-\mathrm{M}\{v\}\in\mathcal{O}_3(x-x_T) \\ v,\ w\in\mathcal{O}_4(x-x_T)\Rightarrow\Delta_{vw}\equiv \mathrm{M}\{v\}-\mathrm{M}\{w\}\in\mathcal{O}_3(x-x_T) \end{array}\right\}\Rightarrow\Delta_{uvw}\equiv\mathrm{M}\{\Delta_{uv}\}-\mathrm{M}\{\Delta_{vw}\}\in\mathcal{O}_2(x-x_T)$$

where $\Delta_{uvw}$ is a quadratic equation one of whose roots is the abscissa $x = x_T$ of the tangent point $T$, inherited from the three quartics that shared it, by means of three overlapping monic polynomial transformations.

Developing the algorithm, after doing heavy calculations derived from the enormous complexity of the polynomials involved, from the overwhelming ocean of numbers and letters emerges the following expression for the quadratic equation $\Delta_{uvw}$ that gives the abscissa of the tangent point $T$: (2.44)

$$a^2 r^4 (a^2 (h^4 \varepsilon^4 (\varepsilon^2-2)+4 h^3 r \varepsilon^3 (2-\varepsilon^2)+h^2 \varepsilon^2 (r^2 (\varepsilon^4+5 \varepsilon^2-12)-k^2 \varepsilon^2 (\varepsilon^2-1) (3 \varepsilon^2-2))+2 h r \varepsilon (2 k^2 \varepsilon^2 (\varepsilon^2-1) (3 \varepsilon^2-2)-r^2 (\varepsilon^4+\varepsilon^2-4))-r^2 (4 k^2 \varepsilon^2 (\varepsilon^2-1) (3 \varepsilon^2-2)+r^2 (2-\varepsilon^4)))-2 a r \varepsilon (b (h^3 \varepsilon^3 (\varepsilon^2-2)+h^2 \varepsilon^2 (k \varepsilon (\varepsilon^2-1) (3 \varepsilon^2-2)+3$$

$$r (2-\varepsilon^2))-3 h r \varepsilon (k \varepsilon (\varepsilon^2-1) (3 \varepsilon^2-2)+r (2-\varepsilon^2))+r^2 (2 k \varepsilon (\varepsilon^2-1) (3 \varepsilon^2-2)+r (2-\varepsilon^2)))+c (h^3 \varepsilon^3 (\varepsilon^2-2)-h^2 \varepsilon^2 (k \varepsilon (\varepsilon^2-1) (3 \varepsilon^2-2)+3 r (\varepsilon^2-2))+3 h r \varepsilon (k \varepsilon (\varepsilon^2-1) (3 \varepsilon^2-2)+r (\varepsilon^2-2))-r^2 (2 k \varepsilon (\varepsilon^2-1) (3 \varepsilon^2-2)+r (\varepsilon^2-2))))-r^2 \varepsilon^2$$

$$(b^2 (3 \varepsilon^4-6 \varepsilon^2+4)+2 b c \varepsilon^2 (4-3 \varepsilon^2)+c^2 (3 \varepsilon^4-6 \varepsilon^2+4)) (h^2 \varepsilon^2-2 h r \varepsilon+r^2)) (h \varepsilon-r)^2 x^2+2 a r^4 (r-h \varepsilon)^3 (a k \varepsilon (\varepsilon^2-1) (2 h^5 \varepsilon^4-12 h^4 r \varepsilon^3-h^3 \varepsilon^2 (2 k^2 \varepsilon^2 (\varepsilon^2-1)+r^2 (\varepsilon^4-2 \varepsilon^2-26))+4 h^2 r \varepsilon (3 k^2 \varepsilon^2 (\varepsilon^2-1)+r^2 (\varepsilon^4-2 \varepsilon^2-6))-h r^2 (24$$

$$k^2 \varepsilon^2 (\varepsilon^2-1)+r^2 (5 \varepsilon^4-10 \varepsilon^2-8))+2 r^3 \varepsilon (8 k^2 (\varepsilon^2-1)+r^2 (\varepsilon^2-2)))+a^2 r (b (h^5 \varepsilon^7-h^4 \varepsilon^5 (4 k (\varepsilon^2-1)+5 r \varepsilon)-h^3 \varepsilon^3 (k^2 \varepsilon^2 (\varepsilon^2-1) (3 \varepsilon^2-4)+20 k r \varepsilon (1-\varepsilon^2)+r^2 (\varepsilon^6-2 \varepsilon^4-7 \varepsilon^2-2))+h^2 r \varepsilon^2 (5 k^2 \varepsilon^2 (\varepsilon^2-1) (3 \varepsilon^2-4)+36 k r \varepsilon (1-\varepsilon^2)+r^2 (3$$

$$\varepsilon^6-6 \varepsilon^4-\varepsilon^2-6))-h r^2 \varepsilon (8 k^2 \varepsilon^2 (\varepsilon^2-1) (3 \varepsilon^2-4)+28 k r \varepsilon (1-\varepsilon^2)+r^2 (3 \varepsilon^6-6 \varepsilon^4+4 \varepsilon^2-6))+r^3 (4 k^2 \varepsilon^2 (\varepsilon^2-1) (3 \varepsilon^2-4)+8 k r \varepsilon (1-\varepsilon^2)+r^2 (\varepsilon^6-2 \varepsilon^4+2 \varepsilon^2-2)))-c (h^5 \varepsilon^7+h^4 \varepsilon^5 (4 k (\varepsilon^2-1)-5 r \varepsilon)-h^3 \varepsilon^3 (k^2 \varepsilon^2 (\varepsilon^2-1) (3 \varepsilon^2-4)+20 k r \varepsilon$$

$$(\varepsilon^2-1)+r^2 (\varepsilon^6-2 \varepsilon^4-7 \varepsilon^2-2))+h^2 r \varepsilon^2 (5 k^2 \varepsilon^2 (\varepsilon^2-1) (3 \varepsilon^2-4)+36 k r \varepsilon (\varepsilon^2-1)+r^2 (3 \varepsilon^6-6 \varepsilon^4-\varepsilon^2-6))-h r^2 \varepsilon (8 k^2 \varepsilon^2 (\varepsilon^2-1) (3 \varepsilon^2-4)+28 k r \varepsilon (\varepsilon^2-1)+r^2 (3 \varepsilon^6-6 \varepsilon^4+4 \varepsilon^2-6))+r^3 (4 k^2 \varepsilon^2 (\varepsilon^2-1) (3 \varepsilon^2-4)+8 k r \varepsilon (\varepsilon^2-1)+r^2 (\varepsilon^6-2 \varepsilon^4+2$$

$$\varepsilon^2-2))))+2 a^2 r^2 \varepsilon^2 (c (h^4 \varepsilon^5+2 h^3 \varepsilon^3 (k (\varepsilon^2-1)-2 r \varepsilon)+2 h^2 r \varepsilon^2 (3 r \varepsilon-4 k (\varepsilon^2-1))+2 h r^2 \varepsilon (5 k (\varepsilon^2-1)-2 r \varepsilon)-r^3 (4 k (\varepsilon^2-1)-r \varepsilon))-b^2 (h^4 \varepsilon^5-2 h^3 \varepsilon^3 (k (\varepsilon^2-1)+2 r \varepsilon)+2 h^2 r \varepsilon^2 (4 k (\varepsilon^2-1)+3 r \varepsilon)-2 h r^2 \varepsilon (5 k (\varepsilon^2-1)+2 r$$

$$\varepsilon)+r^3 (4 k (\varepsilon^2-1)+r \varepsilon)))+r^3 \varepsilon^4 (b-c) (b^2 \varepsilon^2+2 b c (2-\varepsilon^2)+c^2 \varepsilon^2) (h^3 \varepsilon^3-3 h^2 r \varepsilon^2+3 h r^2 \varepsilon-r^3)) x-r^4 (a^4 (h^8 \varepsilon^6-8 h^7 r \varepsilon^5-2 h^6 \varepsilon^4 (k^2 \varepsilon^2 (\varepsilon^2-1)+r^2 (\varepsilon^2-13))+4 h^5 r \varepsilon^3 (4 k^2 \varepsilon^2 (\varepsilon^2-1)+r^2 (3 \varepsilon^2-11))+h^4 \varepsilon^2 (k^4 \varepsilon^4 (\varepsilon^2-1)^2+2 k^2 r^2$$

$$\varepsilon^2 (\varepsilon^2-1) (\varepsilon^4-2 \varepsilon^2-25)+r^4 (\varepsilon^4-32 \varepsilon^2+41))-4 h^3 r \varepsilon (2 k^4 \varepsilon^4 (\varepsilon^2-1)^2+k^2 r^2 \varepsilon^2 (\varepsilon^2-1) (3 \varepsilon^4-6 \varepsilon^2-19)+r^4 (\varepsilon^4-12 \varepsilon^2+5))+2 h^2 r^2 (12 k^4 \varepsilon^4 (\varepsilon^2-1)^2+k^2 r^2 \varepsilon^2 (\varepsilon^2-1) (13 \varepsilon^4-26 \varepsilon^2-28)+r^4 (3 \varepsilon^4-21 \varepsilon^2+2))-4 h r^3 \varepsilon (8 k^4 \varepsilon^2 (\varepsilon^2-1)^2+2 k^2$$

$$r^2 (\varepsilon^2-1) (3 \varepsilon^4-6 \varepsilon^2-2)+r^4 (\varepsilon^4-5))+r^4 (16 k^4 \varepsilon^4 (\varepsilon^2-1)^2+8 k^2 r^2 \varepsilon^2 (\varepsilon^2-1) (\varepsilon^2-2)+r^4 (\varepsilon^4-4)))-2 a^3 r (b (2 h^7 \varepsilon^6-14 h^6 r \varepsilon^5-2 h^5 \varepsilon^4 (k^2 \varepsilon^2 (\varepsilon^2-1)+r^2 (\varepsilon^2-20))+h^4 r \varepsilon^3 (14 k^2 \varepsilon^2 (\varepsilon^2-1)+k r \varepsilon^3 (1-\varepsilon^2) (\varepsilon^2-2)+10 r^2 (\varepsilon^2-6))-h^3 r^2 \varepsilon^2$$

$$(38 k^2 \varepsilon^2 (\varepsilon^2-1)+5 k r \varepsilon^3 (1-\varepsilon^2) (\varepsilon^2-2)+10 r^2 (2 \varepsilon^2-5))+h^2 r^3 \varepsilon (50 k^2 \varepsilon^2 (\varepsilon^2-1)+9 k r \varepsilon^3 (1-\varepsilon^2) (\varepsilon^2-2)+2 r^2 (10 \varepsilon^2-11))-h r^4 (32 k^2 \varepsilon^2 (\varepsilon^2-1)+7 k r \varepsilon^3 (1-\varepsilon^2) (\varepsilon^2-2)+2 r^2 (5 \varepsilon^2-2))+2 r^5 \varepsilon (4 k^2 (\varepsilon^2-1)+k r \varepsilon (1-\varepsilon^2)$$

$$(\varepsilon^2-2)+r^2))+c (2 h^7 \varepsilon^6-14 h^6 r \varepsilon^5-2 h^5 \varepsilon^4 (k^2 \varepsilon^2 (\varepsilon^2-1)+r^2 (\varepsilon^2-20))+h^4 r \varepsilon^3 (14 k^2 \varepsilon^2 (\varepsilon^2-1)+k r \varepsilon^3 (\varepsilon^2-1) (\varepsilon^2-2)+10 r^2 (\varepsilon^2-6))-h^3 r^2 \varepsilon^2 (38 k^2 \varepsilon^2 (\varepsilon^2-1)+5 k r \varepsilon^3 (\varepsilon^2-1) (\varepsilon^2-2)+10 r^2 (2 \varepsilon^2-5))+h^2 r^3 \varepsilon (50 k^2 \varepsilon^2 (\varepsilon^2-1)+9 k r$$

$$\varepsilon^3 (\varepsilon^2-1) (\varepsilon^2-2)+2 r^2 (10 \varepsilon^2-11))-h r^4 (32 k^2 \varepsilon^2 (\varepsilon^2-1)+7 k r \varepsilon^3 (\varepsilon^2-1) (\varepsilon^2-2)+2 r^2 (5 \varepsilon^2-2))+2 r^5 \varepsilon (4 k^2 (\varepsilon^2-1)+k r \varepsilon (\varepsilon^2-1) (\varepsilon^2-2)+r^2)))+a^2 r^2 (b^2 (h^6 \varepsilon^6 (\varepsilon^2+6)-6 h^5 r \varepsilon^5 (\varepsilon^2+6)-h^4 \varepsilon^4 (k^2 \varepsilon^2 (\varepsilon^2-1) (3 \varepsilon^2-2)+r^2 (\varepsilon^4-12$$

$$\varepsilon^2-90))+2 h^3 r \varepsilon^3 (3 k^2 \varepsilon^2 (\varepsilon^2-1) (3 \varepsilon^2-2)+2 r^2 (\varepsilon^4-2 \varepsilon^2-30))-h^2 r^2 \varepsilon^2 (13 k^2 \varepsilon^2 (\varepsilon^2-1) (3 \varepsilon^2-2)+3 r^2 (2 \varepsilon^4+\varepsilon^2-30))+2 h r^3 \varepsilon (6 k^2 \varepsilon^2 (\varepsilon^2-1) (3 \varepsilon^2-2)+r^2 (2 \varepsilon^4+3 \varepsilon^2-18))-r^4 (4 k^2 \varepsilon^2 (\varepsilon^2-1) (3 \varepsilon^2-2)+r^2 (\varepsilon^4+2 \varepsilon^2-6)))-2 b c (h^6$$

$$\varepsilon^6 (\varepsilon^2-6)+6 h^5 r \varepsilon^5 (6-\varepsilon^2)-h^4 \varepsilon^4 (3 k^2 \varepsilon^2 (\varepsilon^2-1) (\varepsilon^2-2)+r^2 (\varepsilon^4-16 \varepsilon^2+86))+2 h^3 r \varepsilon^3 (9 k^2 \varepsilon^2 (\varepsilon^2-1) (\varepsilon^2-2)+2 r^2 (\varepsilon^4-6 \varepsilon^2+26))-3 h^2 r^2 \varepsilon^2 (13 k^2 \varepsilon^2 (\varepsilon^2-1) (\varepsilon^2-2)+r^2 (2 \varepsilon^4-7 \varepsilon^2+22))+2 h r^3 \varepsilon (18 k^2 \varepsilon^2 (\varepsilon^2-1) (\varepsilon^2-2)+r^2 (2 \varepsilon^4-5$$

$$\varepsilon^2+10))-r^4 (12 k^2 \varepsilon^2 (\varepsilon^2-1) (\varepsilon^2-2)+r^2 (\varepsilon^4-2 \varepsilon^2+2)))+c^2 (h^6 \varepsilon^6 (\varepsilon^2+6)-6 h^5 r \varepsilon^5 (\varepsilon^2+6)-h^4 \varepsilon^4 (k^2 \varepsilon^2 (\varepsilon^2-1) (3 \varepsilon^2-2)+r^2 (\varepsilon^4-12 \varepsilon^2-90))+2 h^3 r \varepsilon^3 (3 k^2 \varepsilon^2 (\varepsilon^2-1) (3 \varepsilon^2-2)+2 r^2 (\varepsilon^4-2 \varepsilon^2-30))-h^2 r^2 \varepsilon^2 (13 k^2 \varepsilon^2 (\varepsilon^2-1) (3 \varepsilon^2-2)+3$$

$$r^2 (2 \varepsilon^4+\varepsilon^2-30))+2 h r^3 \varepsilon (6 k^2 \varepsilon^2 (\varepsilon^2-1) (3 \varepsilon^2-2)+r^2 (2 \varepsilon^4+3 \varepsilon^2-18))-r^4 (4 k^2 \varepsilon^2 (\varepsilon^2-1) (3 \varepsilon^2-2)+r^2 (\varepsilon^4+2 \varepsilon^2-6))))-2 a r^3 \varepsilon (b^3 (h^5 \varepsilon^5 (\varepsilon^2+2)+h^4 \varepsilon^4 (k \varepsilon (\varepsilon^2-1) (\varepsilon^2-2)-5 r (\varepsilon^2+2))+5 h^3 r \varepsilon^3 (2 r (\varepsilon^2+2)-k \varepsilon (\varepsilon^2-1) (\varepsilon^2-2))+h^2$$

$$r^2 \varepsilon^2 (9 k \varepsilon (\varepsilon^2-1) (\varepsilon^2-2)-10 r (\varepsilon^2+2))+h r^3 \varepsilon (5 r (\varepsilon^2+2)-7 k \varepsilon (\varepsilon^2-1) (\varepsilon^2-2))+r^4 (2 k \varepsilon (\varepsilon^2-1) (\varepsilon^2-2)-r (\varepsilon^2+2)))-b^2 c (h^5 \varepsilon^5 (\varepsilon^2-6)+h^4 \varepsilon^4 (3 k \varepsilon (\varepsilon^2-1) (\varepsilon^2-2)+5 r (6-\varepsilon^2))-5 h^3 r \varepsilon^3 (3 k \varepsilon (\varepsilon^2-1) (\varepsilon^2-2)+2 r (6-\varepsilon^2))+h^2 r^2 \varepsilon^2$$

$$(27 k \varepsilon (\varepsilon^2-1) (\varepsilon^2-2)+10 r (6-\varepsilon^2))-h r^3 \varepsilon (21 k \varepsilon (\varepsilon^2-1) (\varepsilon^2-2)+5 r (6-\varepsilon^2))+r^4 (6 k \varepsilon (\varepsilon^2-1) (\varepsilon^2-2)+r (6-\varepsilon^2)))-b c^2 (h^5 \varepsilon^5 (\varepsilon^2-6)-h^4 \varepsilon^4 (3 k \varepsilon (\varepsilon^2-1) (\varepsilon^2-2)+5 r (\varepsilon^2-6))+5 h^3 r \varepsilon^3 (3 k \varepsilon (\varepsilon^2-1) (\varepsilon^2-2)+2 r (\varepsilon^2-6))-h^2 r^2 \varepsilon^2 (27$$

$$k \varepsilon (\varepsilon^2-1) (\varepsilon^2-2)+10 r (\varepsilon^2-6))+h r^3 \varepsilon (21 k \varepsilon (\varepsilon^2-1) (\varepsilon^2-2)+5 r (\varepsilon^2-6))-r^4 (6 k \varepsilon (\varepsilon^2-1) (\varepsilon^2-2)+r (\varepsilon^2-6)))+c^3 (h^5 \varepsilon^5 (\varepsilon^2+2)-h^4 \varepsilon^4 (k \varepsilon (\varepsilon^2-1) (\varepsilon^2-2)+5 r (\varepsilon^2+2))+5 h^3 r \varepsilon^3 (k \varepsilon (\varepsilon^2-1) (\varepsilon^2-2)+2 r (\varepsilon^2+2))-h^2 r^2 \varepsilon^2 (9 k \varepsilon (\varepsilon^2-1)$$

$$(\varepsilon^2-2)+10 r (\varepsilon^2+2))+h r^3 \varepsilon (7 k \varepsilon (\varepsilon^2-1) (\varepsilon^2-2)+5 r (\varepsilon^2+2))-r^4 (2 k \varepsilon (\varepsilon^2-1) (\varepsilon^2-2)+r (\varepsilon^2+2))))+2 r^4 \varepsilon^2 (b^4 \varepsilon^2+2 b^3 c (2-\varepsilon^2)+2 b^2 c^2 \varepsilon^2+2 b c^3 (2-\varepsilon^2)+c^4 \varepsilon^2) (h^4 \varepsilon^4-4 h^3 r \varepsilon^3+6 h^2 r^2 \varepsilon^2-4 h r^3 \varepsilon+r^4))=0$$

or in a more compact form, (2.45)

$$(r^2 \varepsilon^2 (\alpha^4 (\varepsilon^2-2)+\alpha^2 \beta^2 (\varepsilon^2+1)+\beta^4)+x_0^2 (\alpha^4 (\varepsilon^2-1)+\beta^4 (3 \varepsilon^2+1))-2 x_0 x_1 (\alpha^4 (\varepsilon^2-1)+\beta^4 (3 \varepsilon^2+1))+x_1^2 (\alpha^4 (\varepsilon^2-1)+\beta^4 (3 \varepsilon^2+1))+y_0^2 (\alpha^2+\beta^2) (\alpha^2 (2$$

$$\varepsilon^2-1)+\beta^2)-2 y_0 y_1 (\alpha^2+\beta^2) (\alpha^2 (2 \varepsilon^2-1)+\beta^2)+y_1^2 (\alpha^2+\beta^2) (\alpha^2 (2 \varepsilon^2-1)+\beta^2)-\alpha^2 \varepsilon^2 (\alpha^4 (\varepsilon^2-1)+\alpha^2 \beta^2 (\varepsilon^2-1)+2 \beta^4)) x^2-2 (r^2 \varepsilon^2 (x_0 (\alpha^2+\beta^2) (\alpha^2 (\varepsilon^2-1)+2$$

$$\beta^2)-x_1 (\alpha^4+\beta^4))+2 x_0^3 \beta^4 \varepsilon^2+x_0^2 x_1 (\alpha^4 (\varepsilon^2-1)+\beta^4 (1-3 \varepsilon^2))+2 x_0 (\alpha^2 \beta^2 \varepsilon^2 ((y_0-y_1)^2-\beta^2)-x_1^2 (\alpha^4 (\varepsilon^2-1)+\beta^4))+x_1 (x_1^2 (\alpha^4 (\varepsilon^2-1)+\beta^4 (\varepsilon^2+1))+y_0^2 (\alpha^4$$

$$(2 \varepsilon^2-1)+\beta^4)-2 y_0 y_1 (\alpha^4 (2 \varepsilon^2-1)+\beta^4)+y_1^2 (\alpha^4 (2 \varepsilon^2-1)+\beta^4)+\alpha^4 \varepsilon^2 (1-\varepsilon^2) (\alpha^2+\beta^2))) x+r^4 \alpha^4 \varepsilon^2+r^2 (x_0^2 (\alpha^2+\beta^2) (\alpha^2 (\varepsilon^2-1)^2+\beta^2 (3 \varepsilon^2-1))+2 x_0 x_1$$

$$(\alpha^4-\beta^4) (\varepsilon^2-1)-x_1^2 (\alpha^4 (2 \varepsilon^2-1)+\beta^4)-2 \alpha^4 \varepsilon^2 ((y_0-y_1)^2+\beta^2))+x_0^4 \beta^4 \varepsilon^2+x_0^2 (x_1^2 (\alpha^4 (\varepsilon^2-1)+\beta^4 (1-3 \varepsilon^2))+2 \alpha^2 \beta^2 \varepsilon^2 ((y_0-y_1)^2-\beta^2))+2 x_0 x_1^3 (1-\varepsilon^2)$$

$$(\alpha^4-\beta^4)+x_1^4 (\alpha^4 (\varepsilon^2-1)+\beta^4)+x_1^2 (y_0^2 (\alpha^4 (2 \varepsilon^2-1)+\beta^4)-2 y_0 y_1 (\alpha^4 (2 \varepsilon^2-1)+\beta^4)+y_1^2 (\alpha^4 (2 \varepsilon^2-1)+\beta^4)+\alpha^4 \varepsilon^2 (1-\varepsilon^2) (\alpha^2+\beta^2))+\alpha^4 \varepsilon^2 (y_0^4-4 y_0^3 y_1+2 y_0^2$$

$$(3 y_1^2-\beta^2)+4 y_0 y_1 (\beta^2-y_1^2)+(y_1^2-\beta^2)^2)=0$$

But the equation (2.27) must have a double root, as the tangency point is the limit to which tend both points of intersection of the circle and the ellipse. Therefore, we can replace the quartic (2.27) by its derivative, which is a cubic: (2.46)

$$u'(x)=4 a^2 \varepsilon^6 (h \varepsilon-r)^4 x^3+12 a \varepsilon^4 (r-h \varepsilon)^3 (a k (\varepsilon^2-1) (h \varepsilon-2 r)+r (b-c) (r-h \varepsilon)) x^2-4 (a^2 (h^4 \varepsilon^4 (\varepsilon^2-2)+4 h^3 r \varepsilon^3 (2-\varepsilon^2)+h^2 \varepsilon^2 (r^2 (\varepsilon^4+5 \varepsilon^2-12)-k^2$$

$$\varepsilon^2 (\varepsilon^2-1) (3 \varepsilon^2-2))+2 h r \varepsilon (2 k^2 \varepsilon^2 (\varepsilon^2-1) (3 \varepsilon^2-2)-r^2 (\varepsilon^4+\varepsilon^2-4))-r^2 (4 k^2 \varepsilon^2 (\varepsilon^2-1) (3 \varepsilon^2-2)+r^2 (2-\varepsilon^4)))+2 a r \varepsilon (r-h \varepsilon) (b (h^2 \varepsilon^2 (\varepsilon^2-2)-2 h \varepsilon$$

$$(k \varepsilon (\varepsilon^2-1)+r (\varepsilon^2-2))+r (4 k \varepsilon (\varepsilon^2-1)+r (\varepsilon^2-2)))+c (h^2 \varepsilon^2 (\varepsilon^2-2)+2 h \varepsilon (k \varepsilon (\varepsilon^2-1)+r (2-\varepsilon^2))-r (4 k \varepsilon (\varepsilon^2-1)+r (2-\varepsilon^2))))-4 r^2 \varepsilon^2 (b^2-b c \varepsilon^2+c^2) (h$$

$$\varepsilon-r)^2) (h \varepsilon-r)^2 x+4 (a^2 k \varepsilon (\varepsilon^2-1) (h \varepsilon-2 r) (h^4 \varepsilon^3-4 h^3 r \varepsilon^2+h^2 \varepsilon (r^2 (\varepsilon^2+5)-k^2 \varepsilon^2 (\varepsilon^2-1))+2 h r (2 k^2 \varepsilon^2 (\varepsilon^2-1)-r^2 (\varepsilon^2+1))-r^2 \varepsilon (4 k^2$$

$$(\varepsilon^2-1)-r^2))+a r (b (h^4 \varepsilon^4-2 h^3 \varepsilon^3 (k \varepsilon (\varepsilon^2-1)+2 r)+h^2 \varepsilon^2 (k^2 \varepsilon^2 (\varepsilon^2-1)+8 k r \varepsilon (\varepsilon^2-1)+r^2 (\varepsilon^2+7))-2 h r \varepsilon (2 k^2 \varepsilon^2 (\varepsilon^2-1)+5 k r \varepsilon (\varepsilon^2-1)+r^2$$

$$(\varepsilon^2+3))+r^2 (4 k^2 \varepsilon^2 (\varepsilon^2-1)+4 k r \varepsilon (\varepsilon^2-1)+r^2 (\varepsilon^2+2)))-c (h^4 \varepsilon^4+2 h^3 \varepsilon^3 (k \varepsilon (\varepsilon^2-1)-2 r)+h^2 \varepsilon^2 (k^2 \varepsilon^2 (\varepsilon^2-1)+8 k r \varepsilon (1-\varepsilon^2)+r^2 (\varepsilon^2+7))-2 h r \varepsilon (2$$

$$k^2 \varepsilon^2 (\varepsilon^2-1)+5 k r \varepsilon (1-\varepsilon^2)+r^2 (\varepsilon^2+3))+r^2 (4 k^2 \varepsilon^2 (\varepsilon^2-1)+4 k r \varepsilon (1-\varepsilon^2)+r^2 (\varepsilon^2+2)))) (h \varepsilon-r)-2 r^2 \varepsilon (b^2 (h \varepsilon-r)^2+2 b c k \varepsilon (1-\varepsilon^2) (h \varepsilon-2 r)-c^2 (h$$

$$\varepsilon-r)^2) (h \varepsilon-r)^2) (h \varepsilon-r)=0$$

Using (2.46) instead of (2.27), the algorithm to obtain the quadratic equation by iteration of the TOP would be: (2.47)

$$\left.\begin{array}{l} u\in\mathscr{O}_4(x-x_T)^2 \Rightarrow u'\in\mathscr{O}_3(x-x_T) \\ v,\ w\in\mathscr{O}_4(x-x_T) \Rightarrow \Delta_{vw}\equiv \mathrm{M}\{v\}-\mathrm{M}\{w\}\in\mathscr{O}_3(x-x_T) \end{array}\right\}\Rightarrow \Delta_{u'vw}\equiv \mathrm{M}\{u'\}-\mathrm{M}\{\Delta_{vw}\}\in\mathscr{O}_2(x-x_T)$$

Thus, we have the abscise $x_T$ as a function of $a$, $b$, $c$, $r$, $\varepsilon$. Replacing $x_T$ in the equation of $\mathscr{C}_1$ it is obtained the corresponding ordinate $y_T$, also as a function of $a$, $b$, $c$, $r$, $\varepsilon$.

## 2.6. Defining the problem.

Let´s gather now a complete set of independent equations to define the problem, translating geometrical meanings to algebraic descriptions.

❶ Condition of **tangency of $\mathscr{C}_2$ and $\mathscr{C}_3$ with each other and with $\mathscr{E}$.** Fusing (1.13) and (1.17):

$$r=\beta\varepsilon,\ \varepsilon\in\left(\frac{1}{\sqrt{2}},1\right) \tag{2.48}$$

❷ Condition of **tangency among the ellipse $\mathscr{E}$ and the sides *a*, *b* and *c* of △*ABC*.** Unifying (2.12) and (2.24):

$$\alpha^2+\beta^2=k^2\left(1-\frac{r}{h\varepsilon-r}\right)^2+\left(h-\frac{r}{\varepsilon}\right)^2 \tag{2.49}$$

where $h$ and $k$ are constants that depend only of $a$, $b$, $c$ according to (2.21) and (2.22).

❸ Condition of **tangency among the outside-circumference $\mathscr{C}_1$ and the ellipse $\mathscr{E}$.**

There are three ways of setting this condition, imposing that: $\mathscr{C}_1$ and $\mathscr{E}$ share the tangent point $T$, or the normal line at $T$ contains the center $C_1$ of $\mathscr{C}_1$, or **$\mathscr{C}_1$ and $\mathscr{E}$ share the tangent line at *T*.** It will be use (2.44) and the equation of $\mathscr{C}_1$.

❹ Condition of **tangency between the outside-circumference $\mathscr{C}_1$ and the sides *b* and *c* of △*ABC*.**

$$x_1=\frac{|c-b|r}{a} \qquad y_1=\frac{(c+b)r}{a} \tag{2.04}$$

❺ Condition of that the distance between the center $C_1$ of $\mathscr{C}_1$ and the tangent point $T$ is $r$.

$$(x_T-x_1)^2+(y_T-y_1)^2=r^2 \tag{2.50}$$

❻ Relationship among the 3 parameters involved in the problem (lengths of the semi-axis of $\mathscr{E}$ and its eccentricity):

$$\varepsilon=\sqrt{1-\frac{\beta^2}{\alpha^2}} \tag{1.11}$$

The condition of equal radius for the three circumferences it is yet included implicitly in the above equations since are used the same variable $r$ for all them.

This is the simplest system of equations for describing the Masayoshi´s problem, but subjacent to its apparent simplicity is hidden an extremely complexity that outcrops when trying to removing the involved parameters, as shown below.

## 2.7. Removing the parameters.

From (1.11) and (2.48),

$$\alpha^2 = \frac{r^2}{\varepsilon^2(1-\varepsilon^2)},\ \varepsilon \in \left(\frac{1}{\sqrt{2}},\ 1\right) \tag{2.51}$$

and from (2.48),

$$\beta = \frac{r}{\varepsilon},\ \varepsilon \in \left(\frac{1}{\sqrt{2}},\ 1\right) \tag{2.52}$$

Now only remains to eliminate the parameter $\varepsilon$, and there is just a way to try it: by means of (2.49), (2.51) and (2.52). Removing $\alpha$ and $\beta$ in (2.12) by means of (2.49) and (2.51), and applying (2.24) to the right side of (2.12),

$$\frac{r^2}{\varepsilon^2(1-\varepsilon^2)} + \left(\frac{r}{\varepsilon}\right)^2 = k^2\left(1-\frac{r}{h\varepsilon - r}\right)^2 + \left(h-\frac{r}{\varepsilon}\right)^2 \tag{2.53}$$

and after a tedious series of operations, expanding and simplifying by means of (2.21), (2.22) and (2.25), it is obtained,

$$\left.\begin{aligned} Q &\equiv \frac{r}{a} \\ R &\equiv \frac{r}{h} \end{aligned}\right\} \quad \varepsilon^6 - 4R\varepsilon^5 + \left\{4\left[Q^2+R^2\right]-1\right\}\varepsilon^4 + 4R\left[1-2Q^2\right]\varepsilon^3 - 4R^2\varepsilon^2 + 4Q^2R^2 = 0 \tag{2.54}$$

Here it is, eureka!: the asymmetric case is solvable iff the sextic (2.54) is solvable algebraically.

There is no further way to eliminate $\varepsilon$, since there are no more independent equations that interrelate these parameters than those that have been used yet. It is, therefore, likely impossible to eliminate the parameter $\varepsilon$ and, consequently, therefore, it is most likely that nor be possible to get a relationship between the variables *a*, *b*, *c* and *r* that doesn´t involve any of the parameters, $\alpha$, $\beta$ or $\varepsilon$.

See appendix C for the analysis of solubility of this equation.

## 2.8. Final relationship between *a*, *b*, *c* and *r*.

If would be possible to solve the sextic (2.54), denoting the right root (that corresponding to the eccentricity of $\mathcal{E}$) by

$$\varepsilon_r \equiv \varepsilon(a, b, c, r) \tag{2.55}$$

the coordinates of *T* would be also functions of *a*, *b*, *c*, *r*, and (2.50) would state an *implicit relationship* among *r* and *a*, *b*, *c*:

$$\left[x_T(r) - x_1(r)\right]^2 + \left[y_T(r) - y_1(r)\right]^2 = r^2 \tag{2.56}$$

where $x_T$ is one of the two roots of (2.44), $y_T$ is the function obtained substituting $x_T$ in the equation of $\mathcal{C}_1$, $x_1$ and $y_1$ the functions defined in (2.04).

## 2.9. Approximate solution by iteration.

Writing (2.56) in the form,

$$r = \sqrt{\left[x_T(r) - x_1(r)\right]^2 + \left[y_T(r) - y_1(r)\right]^2} \equiv \rho(a, b, c, \varepsilon_r, r) \tag{2.57}$$

the algorithm for *n* iterations would be,

$$\begin{cases} r_1 = \rho(r_0) \\ r_2 = \rho(r_1) \\ \quad \vdots \\ r_n = \rho(r_{n-1}) \end{cases} \tag{2.58}$$

where $r_0$ is an *initial value*, arbitrarily chosen, but in the *range of convergence* of the algorithm. The value of *n* depends on the desired accuracy and the algorithm terminates when the values of *r* become stable for a certain number of digits.

Obviously, if is known the value of *r* it is trivial to define one (or any number) relationships between the sides of the triangle and *r*.

**2.9. An alternative way.**

However, by rearranging the equation (2.54) according to the powers of *r*, we get a complete quartic in *r*,

$$\boxed{\sum_{p=0}^{4} a_p \varepsilon^p = 0} \left\{ \begin{array}{l} a_0 \equiv -\left(\dfrac{bc}{2}\right)^2 \varepsilon^4 (1-\varepsilon^2) \\ a_1 \equiv abc\varepsilon^3 (1-\varepsilon^2) \\ a_2 \equiv \left[\left(a^2+h^2\right)\varepsilon^2 - a^2\right]\varepsilon^2 \\ a_3 \equiv -2h\varepsilon^3 \\ a_4 \equiv 1 \end{array} \right. \qquad (2.59)$$

equation which always has algebraic solution. It would, therefore, be necessary to slightly modify the conditions of the statement so that the problem had a solution, simply by changing the variable *r* by the parameter *ε*.

Applying the formula (A.02) given in the appendix A, its roots have the following explicit algebraic expression:

$$\left. \begin{array}{l} \boxed{r = \dfrac{1}{2}\left[ -\dfrac{a_3}{2} + \sqrt{\dfrac{C}{3\sqrt[3]{2}} + E} \pm \sqrt{\dfrac{a_3^2}{2} - \dfrac{4a_2}{3} - \dfrac{C}{3\sqrt[3]{2}} - D + F} \right]} \\ \boxed{r = \dfrac{1}{2}\left[ -\dfrac{a_3}{2} - \sqrt{\dfrac{C}{3\sqrt[3]{2}} + E} \pm \sqrt{\dfrac{a_3^2}{2} - \dfrac{4a_2}{3} - \dfrac{C}{3\sqrt[3]{2}} - D - F} \right]} \end{array} \right\} \left\{ \begin{array}{l} A \equiv 2a_2^3 + 9\left[3(a_1^2 + a_0 a_3^2) - (8a_0 + a_1 a_3)a_2\right] \\ B \equiv a_2^2 + 3(4a_0 - a_1 a_3) \\ C \equiv \sqrt[3]{A + \sqrt{A^2 - 4B^3}} \\ D \equiv \dfrac{\sqrt[3]{2}}{3}\dfrac{B}{C} \\ E \equiv D + \dfrac{a_3^2}{4} - \dfrac{2a_2}{3} \\ F \equiv \dfrac{a_2 a_3 - 2a_1 - \dfrac{a_3^3}{4}}{\sqrt{\dfrac{C}{3\sqrt[3]{2}} + E}} \end{array} \right.$$

Substituting in (A.02) the generic coefficients of the quartic by the specifics ones given by (2.59),

$$\left\{ \begin{array}{l} a_0 \equiv -\left(\dfrac{bc}{2}\right)^2 \varepsilon^4 (1-\varepsilon^2) \\ a_1 \equiv abc\varepsilon^3 (1-\varepsilon^2) \\ a_2 \equiv \left[\left(a^2+h^2\right)\varepsilon^2 - a^2\right]\varepsilon^2 \\ a_3 \equiv -2h\varepsilon^3 \\ a_4 \equiv 1 \end{array} \right\} \Rightarrow \left\{ \begin{array}{l} A \equiv 2a_2^3 + 9\left[3(a_1^2 + a_0 a_3^2) - (8a_0 + a_1 a_3)a_2\right] \\ B \equiv a_2^2 + 3(4a_0 - a_1 a_3) \\ C \equiv \sqrt[3]{A + \sqrt{A^2 - 4B^3}} \\ D \equiv \dfrac{\sqrt[3]{2}}{3}\dfrac{B}{C} \\ E \equiv D + \dfrac{a_3^2}{4} - \dfrac{2a_2}{3} \\ F \equiv \dfrac{a_2 a_3 - 2a_1 - \dfrac{a_3^3}{4}}{\sqrt{\dfrac{C}{3\sqrt[3]{2}} + E}} \end{array} \right. \qquad (2.60)$$

And after some boring simplifications, we obtain as value of the constants $A$, $B$, ..., $F$,

$$\begin{cases} A \equiv \left\{2\left[(a^2+h^2)\varepsilon^2-a^2\right]^3+9b^2c^2(1-\varepsilon^2)\left\{3(a^2-a^2\varepsilon^2-h^2\varepsilon^4)+2(1+\varepsilon^2)\left[(a^2+h^2)\varepsilon^2-a^2\right]\right\}\right\}\varepsilon^6 \\ B \equiv \left\{\left[(a^2+h^2)\varepsilon^2-a^2\right]^2-3b^2c^2(1-\varepsilon^2)(1-2\varepsilon^2)\right\}\varepsilon^4 \\ C \equiv \sqrt[3]{A+\sqrt{A^2-4B^3}} \\ D \equiv \dfrac{\sqrt[3]{2}}{3}\dfrac{B}{C} \\ E \equiv D+h^2\varepsilon^6-\dfrac{2\left[(a^2+h^2)\varepsilon^2-a^2\right]\varepsilon^2}{3} \\ F \equiv \dfrac{2\left\{h^3\varepsilon^6-h\left[(a^2+h^2)\varepsilon^2-a^2\right]\varepsilon^2-abc(1-\varepsilon^2)\right\}\varepsilon^3}{\sqrt{\dfrac{C}{3\sqrt[3]{2}}+E}} \end{cases} \quad (2.61)$$

## 2.10. A numeric example.

Given the following lengths for the sides of $\triangle ABC$,

$$\triangle ABC \begin{cases} a \equiv 2.8939431 \\ b \equiv 1.0591663 \\ c \equiv 2.6931530 \end{cases} \quad (2.62)$$

the values for the remaining parameters involved in the problem, with a 4-decimal precision, are:

$$\begin{cases} \left.\begin{aligned} \alpha &= 1 \\ \beta &= 0.2431 \end{aligned}\right\} \Rightarrow \varepsilon = 0.9700 \\ r = 0.2358 \\ x_0 = 0.7125 \\ y_0 = 0.7425 \\ x_1 = 0.1331 \\ y_1 = 0.3057 \\ x_T = 0.1698 \\ y_T = 0.5386 \end{cases} \quad (2.63)$$

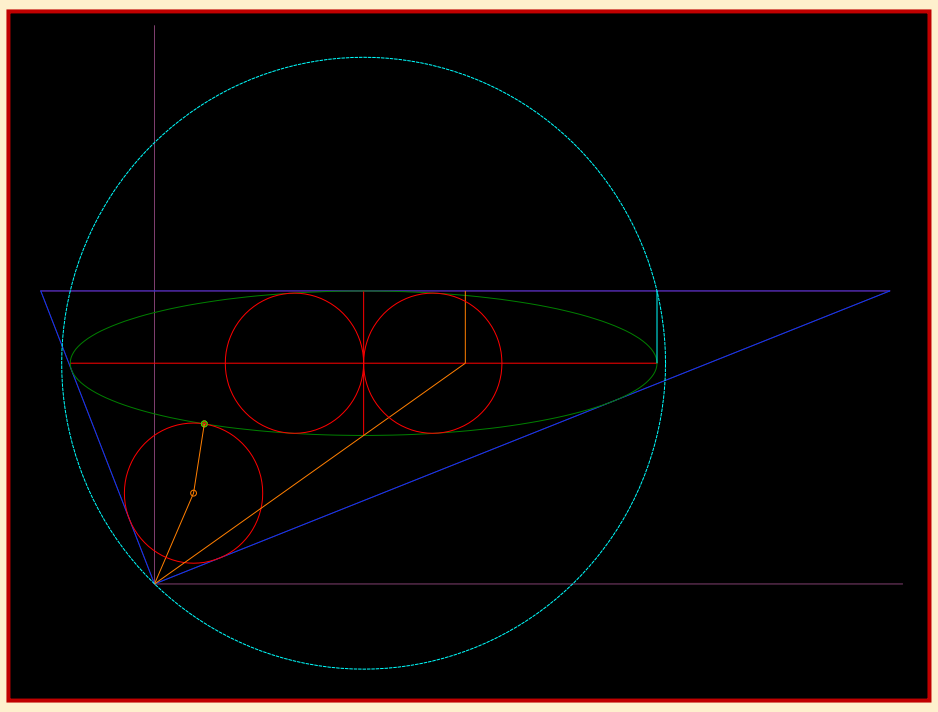

(2.64)

## Appendix A.

**An alternative representation for the explicit relationship given by (1.60) and (1.61).** Expressing (1.22) as,

$$\varepsilon^4 + a_3\varepsilon^3 + a_2\varepsilon^2 + a_1\varepsilon + a_0 = 0 \qquad \text{(A.01)}$$

its roots have the following explicit algebraic expression, (A.02):

$$\left.\begin{array}{l} \varepsilon = \dfrac{1}{2}\left[-\dfrac{a_3}{2} + \sqrt{\dfrac{C}{3\sqrt[3]{2}} + E} \pm \sqrt{\dfrac{a_3^2}{2} - \dfrac{4a_2}{3} - \dfrac{C}{3\sqrt[3]{2}} - D + F}\right] \\ \varepsilon = \dfrac{1}{2}\left[-\dfrac{a_3}{2} - \sqrt{\dfrac{C}{3\sqrt[3]{2}} + E} \pm \sqrt{\dfrac{a_3^2}{2} - \dfrac{4a_2}{3} - \dfrac{C}{3\sqrt[3]{2}} - D - F}\right] \end{array}\right\} \left\{\begin{array}{l} A \equiv 2a_2^3 + 9\left[3(a_1^2 + a_0 a_3^2) - (8a_0 + a_1 a_3)a_2\right] \\ B \equiv a_2^2 + 3(4a_0 - a_1 a_3) \\ C \equiv \sqrt[3]{A + \sqrt{A^2 - 4B^3}} \\ D \equiv \dfrac{\sqrt[3]{2}}{3}\dfrac{B}{C} \\ E \equiv D + \dfrac{a_3^2}{4} - \dfrac{2a_2}{3} \\ F \equiv \dfrac{a_2 a_3 - 2a_1 - \dfrac{a_3^3}{4}}{\sqrt{\dfrac{C}{3\sqrt[3]{2}} + E}} \end{array}\right.$$

Substituting the generic coefficients of the quartic by the specifics ones of our equation, we obtain the value of the constants *A*, *B*, ..., as a function of the constant *p* defined by (1.07):

$$\varepsilon^4 + 2p\varepsilon^3 - \varepsilon^2 - 2p\varepsilon + p^2 = 0 \Rightarrow \left\{\begin{array}{l} a_3 = 2p \\ a_2 = -1 \\ a_1 = -2p \\ a_0 = p^2 \end{array}\right\} \Rightarrow \left\{\begin{array}{l} A \equiv 108p^4 + 144p^2 - 2 \\ B \equiv 24p^2 + 1 \\ C \equiv \sqrt[3]{A + \sqrt{A^2 - 4B^3}} \\ D \equiv \dfrac{\sqrt[3]{2}}{3}\dfrac{B}{C} \\ E \equiv D + p^2 + \dfrac{2}{3} \\ F \equiv \dfrac{2p(1-p^2)}{\sqrt{\dfrac{C}{3\sqrt[3]{2}} + E}} \end{array}\right. \qquad \text{(A.03)}$$

$$p = \sqrt{2} - 1 \Rightarrow \left\{\begin{array}{l} D = \dfrac{\sqrt[3]{2}}{3}\dfrac{B}{C} \left\{\begin{array}{l} B = 73 - 48\sqrt{2} \\ C = \sqrt[3]{2\left[11(103 - 72\sqrt{2}) + 12\sqrt{3(2639 - 1866\sqrt{2})}\right]} \end{array}\right. \\ E = D + 2(2 - \sqrt{2}) - \dfrac{1}{3} \\ F \equiv \dfrac{4(3 - 2\sqrt{2})}{\sqrt{\dfrac{C}{3\sqrt[3]{2}} + E}} \end{array}\right. \qquad \text{(A.04)}$$

Substituting the values given by (1.61) into the explicit expressions for the solutions of the quartic, (1.59), we obtain:

$$\varepsilon = \frac{1}{2}\left[1-\sqrt{2}+\sqrt{\frac{C}{3\sqrt[3]{2}}+E} \pm \sqrt{4(2-\sqrt{2})-\frac{2}{3}-\frac{C}{3\sqrt[3]{2}}-D+F}\right] \Rightarrow \begin{cases} \varepsilon \approx 0.9476620415 \\ \varepsilon \approx 0.1762492425 \end{cases} \quad (A.05)$$

$$\varepsilon = \frac{1}{2}\left[1-\sqrt{2}-\sqrt{\frac{C}{3\sqrt[3]{2}}+E} \pm \sqrt{4(2-\sqrt{2})-\frac{2}{3}-\frac{C}{3\sqrt[3]{2}}-D+F}\right] \Rightarrow \varepsilon \approx -0.9761692044 \pm i\, 0.2726245316 \quad (A.06)$$

and,

$$\frac{b}{r} = 2+\sqrt{2}+\frac{4\sqrt{2}}{1-\sqrt{2}+\sqrt{\frac{C}{3\sqrt[3]{2}}+E}+\sqrt{4(2-\sqrt{2})-\frac{2}{3}-\frac{C}{3\sqrt[3]{2}}-D+F}} \quad (A.07)$$

$$D = \frac{\sqrt[3]{2}}{3}\frac{B}{C} \begin{cases} B = 73-48\sqrt{2} \\ C = \sqrt[3]{2\left[11(103-72\sqrt{2})+12\sqrt{3(2639-1866\sqrt{2})}\right]} \end{cases}$$

$$E = D+2(2-\sqrt{2})-\frac{1}{3} \quad (A.08)$$

$$F \equiv \frac{4(3-2\sqrt{2})}{\sqrt{\frac{C}{3\sqrt[3]{2}}+E}}$$

To shorten the writing of this formula (whose length exceeds the width of this page for a font size readable), have been split into *A*, *B*, ..., *F* parts. This not only represents a greater compactness of the expression but also an economy of calculations because some of these parts is repeated within the expression. Otherwise, writing this steps fully requires too much typing, so it is skipped, but it can be easily verified.

**Appendix B.**

**Developing the equations for the asymmetric case in the reference system used in the symmetric case.**

Let´s consider the same reference system used in the solution of the symmetric case, as shown in the pictures below.

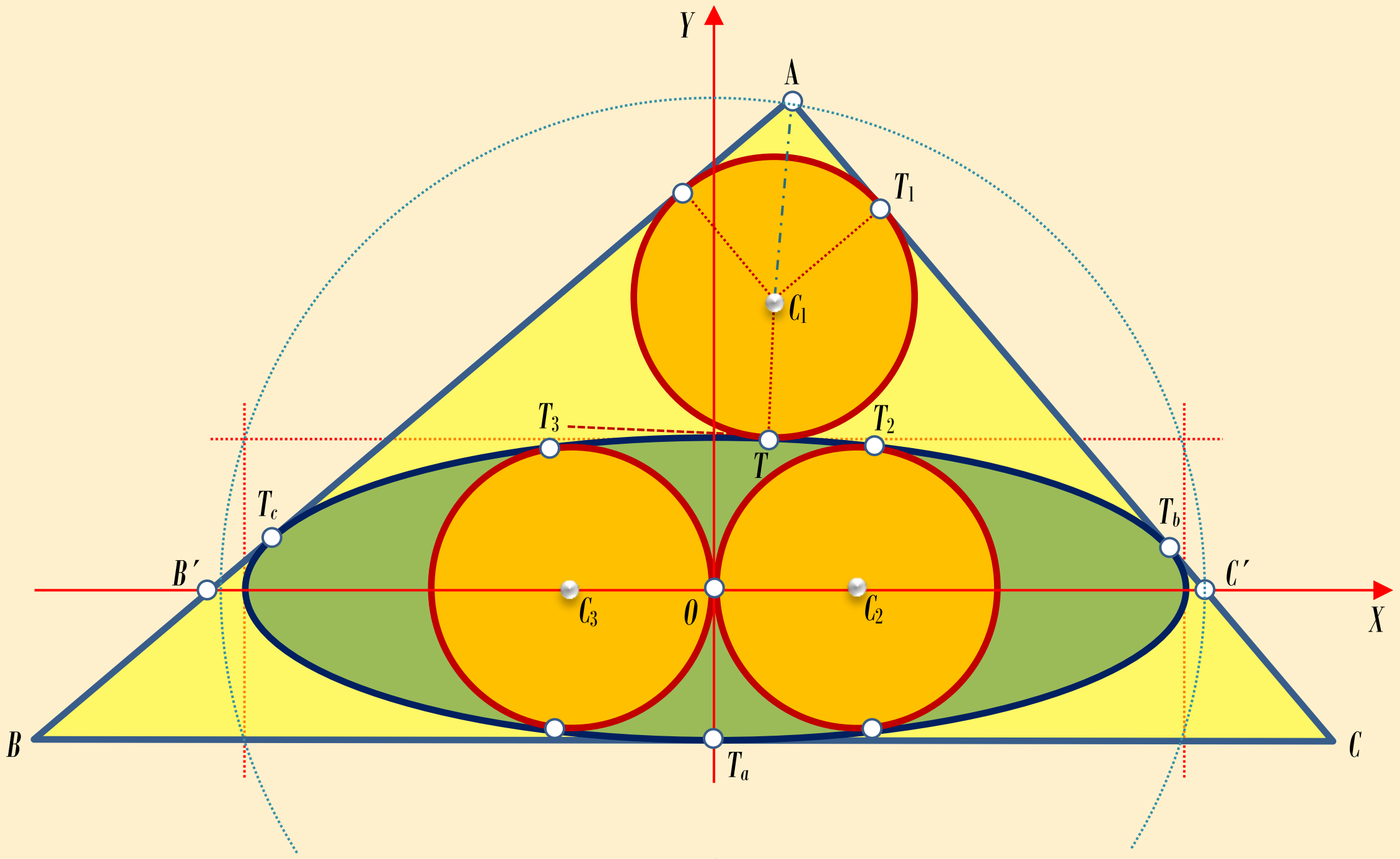

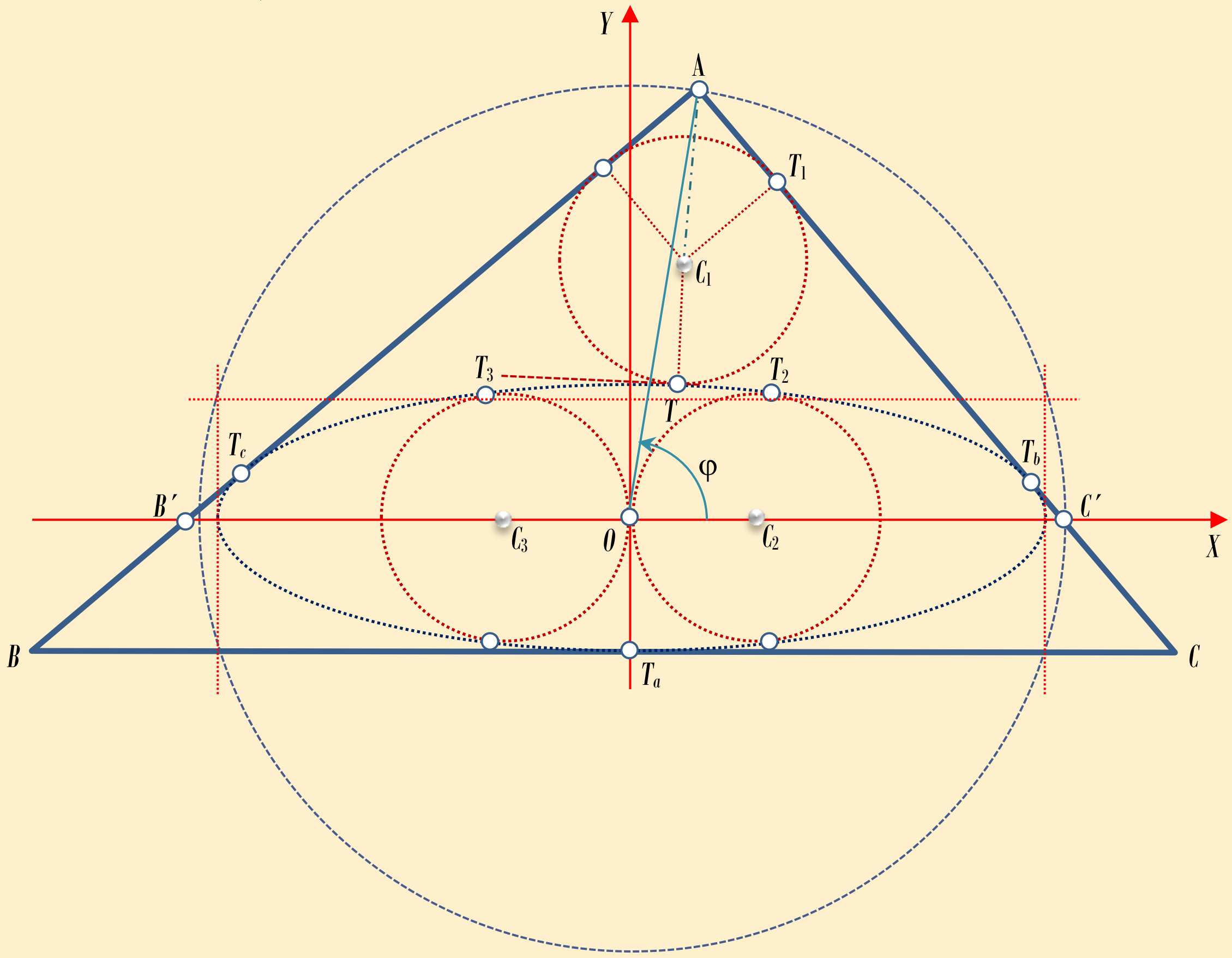

**B.1. Explicit equation of the tangent to the ellipse with slope *m* and ordinate at the origin *n*.**

Let $t \equiv y = mx + n$ a line that intersects the ellipse *e*, (1.01):

$$\left.\begin{array}{l} e \equiv \dfrac{x^2}{\alpha^2} + \dfrac{y^2}{\beta^2} = 1 \\ t \equiv y = mx + n \end{array}\right\} t \cap e \equiv \frac{x^2}{\alpha^2} + \frac{(mx+n)^2}{\beta^2} = 1,$$

i.e.,

$$(\alpha^2 m^2 + \beta^2)x^2 + 2mn\alpha^2 x + \alpha^2(n^2 - \beta^2) = 0 \qquad \text{(B.01)}$$

The *necessary and sufficient condition* for *t* to be tangent to the ellipse is that there is only one point of contact between the line and the curve, which is equivalent to equation (B.01) has only a root (double), ie, its discriminant must be zero:

$$\Delta = (2mn\alpha^2)^2 - 4(\alpha^2 m^2 + \beta^2)\alpha^2(n^2 - \beta^2) = 0 \qquad \text{(B.02)}$$

Operating and simplifying (B.02), we obtain:

$$n = \pm\sqrt{\alpha^2 m^2 + \beta^2} \qquad \text{(B.03)}$$

Thus, the **equation of the tangent of slope *m* to the ellipse** (1.01), for $y > 0$ (positive semi-ellipse), is

$$t \equiv y = mx + \sqrt{\alpha^2 m^2 + \beta^2} \qquad \text{(B.04)}$$

**B.2. Locus of intersection points of pairs of orthogonal tangents: *director circle.***

Let us consider a couple of tangents to the positive semi-ellipse in canonical form, generic and orthogonal each other, (B.04) and

$$s \equiv y = -\frac{1}{m}x + \sqrt{\frac{\alpha^2}{m^2} + \beta^2} \qquad \text{(B.05)}$$

Intersection point of (B.04) and (B.05):

$$s \cap t \equiv \left\{\begin{array}{l} y = -\dfrac{1}{m}x + \sqrt{\dfrac{\alpha^2}{m^2} + \beta^2} \\ y = mx + \sqrt{\alpha^2 m^2 + \beta^2} \end{array}\right\} \Rightarrow x = \frac{\sqrt{\alpha^2 m^{-2} + \beta^2} - \sqrt{\alpha^2 m^2 + \beta^2}}{m + m^{-1}} \qquad \text{(B.06)}$$

Let $A(x_A, y_A) \equiv s \cap t.$

Multiplying by *m* the numerator and denominator of (B.06), we have:

$$x_A = \frac{\sqrt{\alpha^2 + \beta^2 m^2} - m\sqrt{\alpha^2 m^2 + \beta^2}}{1 + m^2} \qquad \text{(B.07)}$$

$$\because A(x_A, y_A) \in t \equiv y = mx + n,\; y_A = m\frac{\sqrt{\alpha^2+\beta^2 m^2} - m\sqrt{\alpha^2 m^2+\beta^2}}{1+m^2} + \sqrt{\alpha^2 m^2+\beta^2},\text{ i.e.,}$$

$$\boxed{y_A = \frac{m\sqrt{\alpha^2+\beta^2 m^2} + \sqrt{\alpha^2 m^2+\beta^2}}{1+m^2}} \qquad \text{(B.08)}$$

Let us compact the notation through the variable changes, $u \equiv \sqrt{\alpha^2+\beta^2 m^2}$, $v \equiv \sqrt{\alpha^2 m^2+\beta^2}$, and calculate the sum of the squares of the coordinates of *A*, given by (B.07) and (B.08):

$$x_A^2 + y_A^2 = \left(\frac{u-mv}{1+m^2}\right)^2 + \left(\frac{mu+v}{1+m^2}\right)^2 = \frac{\cancel{(1+m^2)}(u^2+v^2)}{(1+m^2)^{\cancel{2}}} = \frac{u^2+v^2}{1+m^2},$$

and undoing the variable changes,

$$x_A^2 + y_A^2 = \frac{u^2+v^2}{1+m^2} = \frac{\left(\sqrt{\alpha^2+\beta^2 m^2}\right)^2 + \left(\sqrt{\alpha^2 m^2+\beta^2}\right)^2}{1+m^2} = \frac{\cancel{(1+m^2)}(\alpha^2+\beta^2)}{\cancel{1+m^2}},$$

i.e.,

$$\boxed{x_A^2 + y_A^2 = \alpha^2+\beta^2} \qquad \text{(B.09)}$$

*The locus of the points from which can be drawn pairs of orthogonal tangents to the canonical ellipse is the circle centered at the origin of coordinates and radius* $R = \sqrt{\alpha^2+\beta^2}$. This circle is called the ***director circle***.

### B.3. Coordinates of the tangent point $T_C$ of the line with slope *m* and the ellipse.

Solving the equation (B.01), under the tangency condition established in (B.02) and the value of *n* according to (B.04), we get the abscissa $x_T$ of the tangent point $T_C$:

$$\left.\begin{array}{l} (\alpha^2 m^2+\beta^2)x^2 + 2mn\alpha^2 x + \alpha^2(n^2-\beta^2) = 0 \\ \Delta = 0 \\ n = \sqrt{\alpha^2 m^2+\beta^2} \end{array}\right\} \Rightarrow x_T = \frac{-2\alpha^2 m\sqrt{\alpha^2 m^2+\beta^2}}{2(\alpha^2 m^2+\beta^2)},$$

or, $x_T = \dfrac{-\alpha^2 m}{\sqrt{\alpha^2 m^2+\beta^2}}$.

By replacing $x_T$ into $t \equiv y = mx+n$, we obtain, $y_T = \dfrac{-\alpha^2 m^2}{\sqrt{\alpha^2 m^2+\beta^2}} + \sqrt{\alpha^2 m^2+\beta^2} = \dfrac{\beta^2}{\sqrt{\alpha^2 m^2+\beta^2}}$.

Gathering the above results, we have the following expression for the **coordinates of $T_C$:**

$$\boxed{\mathrm{T}\left(-\frac{\alpha^2 m}{\sqrt{\alpha^2 m^2+\beta^2}}, \frac{\beta^2}{\sqrt{\alpha^2 m^2+\beta^2}}\right)} \qquad \text{(B.10)}$$

And from (1.11), as a function of the eccentricity of the ellipse,

$$\boxed{x_{\mathrm{T}} = -\frac{m}{1-\varepsilon^2}\cdot\frac{r}{\varepsilon\sqrt{1+\dfrac{m^2}{1-\varepsilon^2}}},\; y_{\mathrm{T}} = \frac{r}{\varepsilon\sqrt{1+\dfrac{m^2}{1-\varepsilon^2}}}} \qquad \text{(B.11)}$$

## B.4. Equation of the bisector of the right angle A.

$$\left.\begin{array}{l} c \equiv y = m(x - x_A) + y_A \quad \Leftrightarrow \dfrac{mx - y - mx_A + y_A}{\sqrt{1+m^2}} = 0 \\ b \equiv y = -\dfrac{1}{m}(x - x_A) + y_A \Leftrightarrow \dfrac{x + my - x_A - my_A}{\sqrt{1+m^2}} = 0 \end{array}\right\} \Leftrightarrow mx - y - mx_A + y_A = \pm(x + my - x_A - my_A)$$

$$mx - y - mx_A + y_A = +(x + my - x_A - my_A) \Leftrightarrow \boxed{(1-m)x + (1+m)y - (1-m)x_A - (1+m)y_A = 0}$$

$$mx - y - mx_A + y_A = -(x + my - x_A - my_A) \Leftrightarrow \boxed{(1+m)x - (1-m)y - (1+m)x_A + (1-m)y_A = 0}$$

For $A$ in the first quadrant its bisector has positive slope, so its implicit equation is,

$$\boxed{(1+m)x - (1-m)y - (1+m)x_A + (1-m)y_A = 0}$$

and explicitly,

$$\boxed{y = \frac{1+m}{1-m}x + y_A - \frac{1+m}{1-m}x_A} \qquad \text{(B.12)}$$

### Another way to deduce the equation of the bisector.

The slope of the bisector is

$$\mu = tg\left(B + \frac{\pi}{4}\right) = \frac{tg\,B + tg\,\frac{\pi}{4}}{1 - tg\,B \cdot tg\,\frac{\pi}{4}} = \frac{1+m}{1-m}.$$

Hence, the equation of the bisector is:

$$y = \mu(x - x_A) + y_A = \frac{1+m}{1-m}x - \frac{1+m}{1-m}x_A + y_A.$$

Substituting $x_A$ by the expression (B.07) and $y_A$ by the right side of (B.04), taking into account that $A \in t$,

$$y = \frac{1+m}{1-m}x - \frac{1+m}{1-m}\cdot x_A + mx_A + n = \frac{1+m}{1-m}x + \left(m - \frac{1+m}{1-m}\right)\cdot x_A + \sqrt{\alpha^2 m^2 + \beta^2} =$$

$$= \frac{1+m}{1-m}x - \frac{\cancel{1+m^2}}{1-m}\frac{\sqrt{\alpha^2 + \beta^2 m^2} - m\sqrt{\alpha^2 m^2 + \beta^2}}{\cancel{1+m^2}} + \sqrt{\alpha^2 m^2 + \beta^2},$$

i.e.,

$$\boxed{y = \frac{1+m}{1-m}x + \frac{\sqrt{\alpha^2 m^2 + \beta^2} - \sqrt{\alpha^2 + \beta^2 m^2}}{1-m}} \qquad \text{(B.13)}$$

or,

$$\boxed{y = \frac{1+m}{1-m}x + \frac{r}{1-m}\cdot\frac{\sqrt{1+m^2-\varepsilon^2} - \sqrt{1+m^2-m^2\varepsilon^2}}{\varepsilon\sqrt{1-\varepsilon^2}}} \qquad \text{(B.14)}$$

**B.5. Coordinates of $C_1$.**

$C_1$ is determined by the intersection of the bisector of $A$, (2.12), and the circumference with center $A$ and radius $r\sqrt{2}$:

$$\left.\begin{array}{l}(x-x_A)^2+(y-y_A)^2=2r^2 \Rightarrow \boxed{y_1=y_A+\sqrt{2r^2-(x_1-x_A)^2}} \\ \boxed{y_1=\frac{1+m}{1-m}x_1+y_A-\frac{1+m}{1-m}x_A}\end{array}\right\} \Leftrightarrow \sqrt{2r^2-(x_1-x_A)^2}=\frac{1+m}{1-m}x_1-\frac{1+m}{1-m}x_A$$

or,

$$2r^2-(x_1-x_A)^2=\left(\frac{1+m}{1-m}\right)^2(x_1-x_A)^2,$$

hence,

$$2r^2=\left[\left(\frac{1+m}{1-m}\right)^2+1\right](x_1-x_A)^2 \Leftrightarrow r^2=\frac{1+m^2}{(1-m)^2}(x_1-x_A)^2,\ if\ \ m>0$$

from which,

$$\boxed{x_1=x_A-\frac{1-m}{\sqrt{1+m^2}}r} \qquad \text{(B.15)}$$

and substituting this expression for $x_1$ into (B.12),

$$\boxed{y_1=y_A-\frac{1+m}{\sqrt{1+m^2}}r} \qquad \text{(B.16)}$$

Replacing in (B.15) and (2.16) the coordinates of $A$ given by (B.07) and (B.08), respectively,

$$(x_1,\ y_1)=\left(\frac{\sqrt{\alpha^2+\beta^2m^2}-m\sqrt{\alpha^2m^2+\beta^2}}{1+m^2}-\frac{1-m}{\sqrt{1+m^2}}r,\ \frac{m\sqrt{\alpha^2+\beta^2m^2}+\sqrt{\alpha^2m^2+\beta^2}}{1+m^2}-\frac{1+m}{\sqrt{1+m^2}}r\right)$$

$$\boxed{\begin{array}{l} x_1=\dfrac{\beta\left(\sqrt{\dfrac{\alpha^2}{\beta^2}+m^2}-m\sqrt{\dfrac{\alpha^2}{\beta^2}m^2+1}\right)}{1+m^2}-\dfrac{1-m}{\sqrt{1+m^2}}r \\ y_1=\dfrac{\beta\left(m\sqrt{\dfrac{\alpha^2}{\beta^2}+m^2}+\sqrt{\dfrac{\alpha^2}{\beta^2}m^2+1}\right)}{1+m^2}-\dfrac{1+m}{\sqrt{1+m^2}}r \end{array}} \qquad \text{(B.17)}$$

or, as a function of the eccentricity of the ellipse:

$$\boxed{\begin{array}{l} x_1=\dfrac{\dfrac{r}{\varepsilon}\left(\sqrt{\dfrac{1}{1-\varepsilon^2}+m^2}-m\sqrt{\dfrac{m^2}{1-\varepsilon^2}+1}\right)-r(1-m)\sqrt{1+m^2}}{1+m^2} \\ y_1=\dfrac{\dfrac{r}{\varepsilon}\left(m\sqrt{\dfrac{1}{1-\varepsilon^2}+m^2}+\sqrt{\dfrac{m^2}{1-\varepsilon^2}+1}\right)-r(1+m)\sqrt{1+m^2}}{1+m^2} \end{array}} \qquad \text{(B.18)}$$

**Appendix C. Analysis of the solubility of the sextic (2.54).**

As is well known, [3] and [4]. by means of a *Tschirnhausen transformation* the general sextic,

$$x^6 + a_5 x^5 + a_4 x^4 + a_3 x^3 + a_2 x^2 + a_1 x + a_0 = 0 \qquad \text{(C.01)}$$

can be reduced to the form,

$$z^6 + z^2 + a_2 z + a_1 = 0 \qquad \text{(C.02)}$$

Thus, via *Tschirnhausen transformation*, the sextic (2.54) can be reduced to the form,

$$z\left[z^5 + z - \left(\frac{2r}{b}\right)^2 - \left(\frac{2r}{c}\right)^2\right] = 0;\ \ b, c, r \in \mathbb{R}^+;\ \ 2r < b < c \qquad \text{(C.03)}$$

All equations of degree $\leq 4$ are algebraically solvable over the complex number field $\mathbb{C}$, but from the work of *Ruffini*, *Abel* and *Galois* it is now well-known that all quintic equations do not [2].

In [1] are stated the *necessary and sufficient conditions* for the quintic factor of (C.03), a *Bring-Jerrard* form or *principal quintic* for the general quintic, be solvable by radicals: the existence of three *rational* numbers, $\epsilon = \pm 1$, $p > 0$ and $q \neq 0$, such that, applied to our specific values,

$$\frac{5q^4(3-4\epsilon\, p)}{p^2+1} = 1,\ \ \frac{4q^5(11\epsilon + 2p)}{p^2+1} = \left(\frac{2r}{b}\right)^2 + \left(\frac{2r}{c}\right)^2 \qquad \text{(C.04)}$$

from where,

$$\frac{4q(11\epsilon + 2p)}{5(3-4\epsilon\, p)} = \left(\frac{2r}{b}\right)^2 + \left(\frac{2r}{c}\right)^2 \qquad \text{(C.05)}$$

that can be written as,

$$\frac{4q(11\epsilon + 2p)}{5(3-4\epsilon\, p)} = 4\left(\frac{a}{bc}\right)^2 r^2 \qquad \text{(C.06)}$$

or,

$$\frac{4q(11\epsilon + 2p)}{5(3-4\epsilon\, p)} = \left(\frac{2r}{h}\right)^2 \qquad \text{(C.07)}$$

where $r$ is the radius common to the three circumferences and $h$ the height relative to the hypotenuse of $\triangle ABC$.

In [4] are given explicit formulas to get the roots for this irreducible particular form of quintic when it is algebraically solvable:

$$z_j = \sum_{k=1}^{4} e^{jk\left(\frac{2\pi i}{5}\right)} u_k,\ \ j = [0, 4] \cap \mathbb{Z} \qquad \text{(C.08)}$$

with,

$$u_1 \equiv \sqrt[5]{\frac{v_1^2 v_3}{\Delta^2}},\ u_2 \equiv \sqrt[5]{\frac{v_3^2 v_4}{\Delta^2}},\ u_3 \equiv \sqrt[5]{\frac{v_2^2 v_1}{\Delta^2}},\ u_4 \equiv \sqrt[5]{\frac{v_4^2 v_2}{\Delta^2}},\ \Delta \equiv p^2+1 \quad \text{(C.09)}$$

$$v_1 \equiv \sqrt{\Delta}+\sqrt{\Delta-\epsilon\sqrt{\Delta}},\ v_2 \equiv -\sqrt{\Delta}-\sqrt{\Delta+\epsilon\sqrt{\Delta}},\ v_3 \equiv -\sqrt{\Delta}+\sqrt{\Delta+\epsilon\sqrt{\Delta}},\ v_4 \equiv \sqrt{\Delta}-\sqrt{\Delta-\epsilon\sqrt{\Delta}} \quad \text{(C.10)}$$

Furthermore, any *irreducible* equation with rational coefficients is solvable by radicals if and only if its *Galois group G* is a subgroup of one of the two transitive subgroups of $S_6$: $J = S_2 \# S_3$ and $K = S_3 \# Z_2$, denoting the *wreath product* of two groups by #, [5].

At first glance, the lack of the first-degree term in (2.54) might seem to facilitate the analysis of solvability of the equation, but is just the opposite because the *resolvent equation* of (C.02) is of the form [1],

$$z^{15}-6a_1^2z^{13}-(42a_1+3)a_1^3z^{12}+\cdots+(144a_1a_2^2-32a_2^4-3)a_1^{15}a_2^2=0 \quad \text{(C.11)}$$

and (C.02) is solvable by radicals and its *Galois group* is a subgroup of *J* iff (G.11) has a rational root, but it is straightforward that this happens if the constant term vanishes because then $z = 0$, i.e., if

$$a_1 = \frac{32a_2^4+3}{144a_2^2} \quad \text{(C.12)}$$

and taking, for instance, $a_2 = ½$, we find that the *Galois group* of the irreducible polynomial $36z^6 + 36z^2 + 18z + 5$ is a subgroup of *J* and solvable. But in our case, $a_1 = 0$, and so (C.02) is *reducible*, therefore this theorem is not applicable.

Moreover, the second condition of (C.04) involves a relationship between *b*, *c* and *r*, which is the final goal of the problem, so this condition neither is implementable.

But, according to *Dummit's analysis* [3], solvable irreducible quintic equations are rare, hence, it is reasonable to expect that solvable irreducible sextics are also rare.

Through the exhaustive analysis carried out we have seen that the insurmountable barrier to solving the asymmetric case is the sextic in $\varepsilon$ (2.54), the only way to eliminate the parameter $\varepsilon$.

## Appendix D. Warnings about the polynomial 6th degree equation (2.54).

One of the roots of the sextic (2.54), expressed in its most compact form,

$$\boxed{\sum_{q=0}^{6} c_q \varepsilon^q = 0} \begin{cases} c_0 \equiv \left(\dfrac{2r^2}{bc}\right)^2 \\ c_1 \equiv 0 \\ c_2 \equiv -\left(\dfrac{2ar}{bc}\right)^2 \\ c_3 \equiv \dfrac{4r(a^2-2r^2)}{abc} \\ c_4 \equiv 4r^2\left[\left(\dfrac{a}{bc}\right)^2 + \dfrac{1}{a^2}\right] - 1 \\ c_5 \equiv -\dfrac{4ar}{bc} \\ c_6 \equiv 1 \end{cases} \qquad \text{(D.01)}$$

is the eccentricity $\varepsilon$ of the ellipse, in terms of $a$, $b$, c and $r$, but this equation could be unsolvable by radicals.

As a "warning to sailors", it can be useful to reproduce the following chain of thoughts: we know that (D.01) has to have, at least, one real solution, and taking into account that not real roots appear only in pairs of conjugate complex, necessarily, the number of real roots have to be even, and by the *uniqueness* of the solution discussed before, we´d can think (wrongly) that the unique real solution of this equation has to have multiplicity 2, 4 or 6, and, therefore, the derivative of this polynomial should preserve the same real root, allowing the substitution of the before equation by the following:

$$\sum_{k=1}^{6} k c_k \varepsilon^{k-1} = 0 \qquad \text{(D.02)}$$

But $c_1 = 0$, therefore, (D.02) can be written as,

$$\sum_{k=2}^{6} k c_k \varepsilon^{k-1} = 0 \qquad \text{(F.03)}$$

At this point, we would think that the lack of the first grade term in $\varepsilon$ in the equation (D.01) should be the key for the door to the solution, since it allow to get down one grade more by elimination of the common factor $\varepsilon$ in (D.03) (it has been proved yet that the solution $\varepsilon = 0$ it isn´t admissible), getting finally the following quartic:

$$\sum_{k=2}^{6} k c_k \varepsilon^{k-2} = \sum_{p=0}^{4} (p+2) c_{p+2} \varepsilon^p = 0 \qquad \text{(D.04)}$$

Unfortunately, the hypothesis that the equation (F.01) has one unique manifold real root can´t be proved, hence, it isn´t warranted to reach the solution trough this path. Indeed, for the numerical example presented in 2.11, the six roots of (2.54) are all real and different.

## Appendix E.

### Exploring another ways.

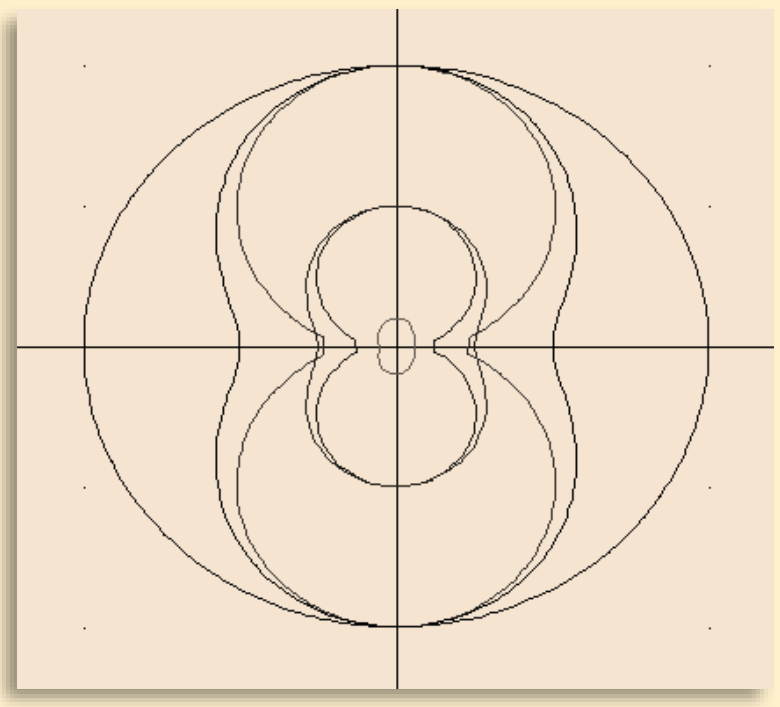

The transformation by *inversion*, the most widely-used method in Sangaku challenges, applied to the ellipse, transform its canonical equation into a quartic equation which is much more complicated.

For instance, the curves on the right picture arise from the mapping

$$(x,\ y) \to \frac{(x, -y)}{x^2 + y^2}$$

by *inversion on the unit circle*, that transforms the ellipse

$$\frac{x^2}{\alpha^2} + \frac{y^2}{\beta^2} = 1$$

into the *quartic*

$$\frac{x^2}{\alpha^2} + \frac{y^2}{\beta^2} = (x^2 + y^2)^2.$$

It is neither possible to simplify the problem trough affine transformations, since an *affinity* that transform the ellipse into a circumference, it also will transform, at the same time, the circumferences into ellipses.

**References.**

[1] Solvable Sextic Equations. C. Boswell and M.L. Glasser. arXiv:math-ph/0504001v1.

[2] R. Bruce King, Beyond the Quartic Equation [Birkhauser, Boston 1996].

[3] D.S. Dummit, Math. Comp 195, 387-401 (1991).

[4] B.K. Spearman and K.S. Williams, Amer. Math. Monthly 101, 986-992 (1994).

[5] J.D. Dixon and B. Mortimer, Permutation Groups,[Springer, New York 1991].